\documentclass[12pt,a4paper,twoside]{article}
\usepackage{amsfonts}
\usepackage{amssymb}
\usepackage[mathcal]{euscript}
\usepackage{ergomath}
\usepackage{comments}
\KillComments
\newKibitzer\DL{D}{red}
\newKibitzer\VM{V}{blue}

\KillComments
\usepackage{smart}
\usepackage[fleqn]{amsmath}
\Aiv
\theorems
\link{equation}{section}
\iffalse
\else
\plainpage
\fi
\margins*[2cm]

\newcommand{\rmname}[1]
  {\expandafter\newcommand \csname #1\endcsname {{\operatorname{#1}}}}
\rmname{Hom}
\rmname{Id}
\rmname{End}
\rmname{ad}
\mathdef \fg     ${\mathfrak{g}}$
\mathdef \fh     ${\mathfrak{h}}$
\let\cal\mathcal
\mathdef \cD    ${\cal D}$
\mathdef \cP    ${\cal P}$
\mathdef\K$\Bbbk$
\mathdef\N$\Bbb N$
\mathdef\R$\Bbb R$
\mathdef\C$\Bbb C$
\mathdef\Z$\Bbb Z$
\mathdef\Zplus$\Bbb Z_+$
\mathdef\Zz$\Bbb Z_2$
\mathdef\La$\Lambda$
\mathdef\One$\mbox{{\rm1\kern-0.37em1}}$
\warndef\harr{\hookrightarrow}
\warndef\arr{\longrightarrow}
\mathdef\sar#1,#2,#3;$#1\stackrel{#2}{\arr}#3$
\mathdef\diff#1${\rm Diff}(#1)$
\mathdef\hom#1#2#3${\rm Hom}_{#1}(#2,#3)$
\mathdef\homk#1#2$\hom
{}{#1}{#2}$
\mathdef\L#1#2$\cal L(#1,#2)$
\mathdef\fam#1#2$\{#1\}_{#2}$
\mathdef\fami#1$\{#1\}_{i\in I}$
\mathdef\famj#1$\{#1\}_{j\in j}$
\warndef\hhat{\widehat}
\let\wh\hhat
\warndef\comp{\mbox{\scriptsize$\circ$}}
\warndef\mult{\centerdot}
\warndef\cdc{,\dots,}
\warndef\la{\langle}
\warndef\ra{\rangle}
\warndef\eop{\hbox{\vrule width 6pt height 6pt depth 0pt}}
\warndef\hotimes{\hhat\otimes}
\expandafter\ifx\csname endproof\endcsname\relax
\newenvironment{proof}{\noindent{\bf Proof}}{\eop\par\medskip}
\fi

\mathdef\0$\bar0$
\mathdef\1$\bar1$
\mathdef \fU     $\mathfrak{U}$
\mathdef\eps$\varepsilon$
\mathdef\I#1$I_\fh^\fg(#1)$
\mathdef\T#1$T_\fh^\fg(#1)$
\mathdef\U#1$U(#1)$
\def\ParagraphSign{\mathhexbox278}
\mathdef\S#1$S(#1)$
\mathdef\cd#1#2#3$#1_{#2}\cdc#1_{#3}$
\mathdef\cu#1#2#3$#1^{#2}\cdc#1^{#3}$
\mathdef\df$\stackrel{{\rm def}}{=}$
\let\la\langle
\let\ra\rangle
\mathdef\kspan#1$\K\la\,#1\,\ra$
\mathdef\Sstar#1$S(#1)^*$
\mathdef\LVect$\text{\rm\bfseries LVect}$
\mathdef\CGLVect$\text{\rm\bfseries CGVect}$
\mathdef\PCGLVect$\text{\rm\bfseries PCGVect}$
\mathdef\bil#1,#2,#3;$\text{\rm Bil}(#1\times#2,#3)$
\mathdef\botimes$\hbox{\rm$\kern.1em\boldsymbol\otimes\kern-.95em\boldsymbol\otimes\kern.01em$}$
\mathdef\hbotimes$\hhat\botimes$
\mathdef\otimesp$\otimes_{\rm p}$
\mathdef\hotimesp$\hotimes_{\rm p}$
\mathdef\diffinf#1${\rm Diff}_\infty(#1)$
\mathdef\A$G_-$
\mathdef\B$H$
\mathdef\chA$\check G_-$
\mathdef\chB$\check H$

\def\M{g}
\begin{document}
\subjclass{17B15 (Primary), 17B65 (Secondary)}
\keywords{Lie superalgebras, representations}
\address{MPIMiS, Inselstr. 22, DE-04103 Leipzig, Germany \\
on leave from Institute for Nuclear Research and Nuclear Energy\\
Tsarigradsko chosse blvd. 72, BG-1784 Sofia, Bulgaria\\
e-mail: {\sl vmolot@inrne.bas.bg}}
\author{Vladimir Molotkov
\thanks{This work was written with financial support of
MPIMiS (Leipzig) and thanks to creative atmosphere of this
Institution. I am especially thankful to A.~Lebedev,  D.~Leites
and Ch.~Sachse for helpful discussion, triggering the idea of
representing $\protect\T g$  and $\protect\I g$
(actions of an element $g\in\fg$ in coinduced and induced \fg-module)
as ``path
integral'' in the action graph, which inspired me to write a
program making all calculations: PCs are very capable playing with
graphs.}
}%
\title{
Explicit Realization of Induced and Coinduced modules 
over Lie Superalgebras by Differential Operators}
\date{September 9, 2005}
\maketitle

This text is an extended version of a part of my 1980 Trieste
preprint \cite{1}. In this preprint I gave, in particular, an
expression for the action of a Lie superalgebra \fg on the
\fg-module coinduced from an \fh-module $V$ for the case where
both \fg and $V$ are finite-dimensional and \fg decomposes (as a
{\it linear space}, not Lie superalgebra) into a direct sum
$\fg_-\oplus\fh$ with $\fg_-$ a Lie superalgebra.
\footnote{As far as I know, the first paper in which there
appeared an explicit expression for the action of a
finite-dimensional Lie algebra $\fg=\fg_-\oplus\fh$ (represented
as a direct sum of Lie subalgebras) by differential operators on
the space $\fg_-^*$ was  I.~Kantor's paper~\cite{2}  of whose
existence I
became aware only last year. In  \cite{3}, I.~Kantor gave an
expression for this action in a more general situation, when (1)
$\fg_-$ may be not a sub{\sl algebra} but only a
finite-dimensional sub{\sl space} of \fg and (2) the Lie algebra
\fh may be infinite-dimensional.

In both \cite{2} and \cite{3}, I.~Kantor only considered Lie algebras and their quasiregular
representations, which are, in terms of coinduced modules, \fg-modules coinduced from
a trivial 1-dimensional 
\fh-module.}

The importance of coinduced and induced \fg-modules is a
consequence of the fact that the vast majority of representations
of Lie superalgebras occurring both in mathematical and physical
practice, belongs to one of these two classes (e.g., the
``universal'' \fg-module $U(\fg)$ and its dual ``couniversal''
module $U(\fg)^*$,
Verma and Harish-Chandra modules, etc.). 
Nevertheless, even some mathematicians working with these
representations do not realize that they are ``speaking prose'',
to say nothing of physicists who study ``representations by
creation and annihilation operators'' while speaking about Verma
modules or their duals.

The standard method used by both mathematicians and physicists to calculate
the action for basis elements of $\fg_-$ in coinduced modules
(where they act as differential operators) in case where \fg is a finite-dimensional Lie algebra,
was to calculate the differential of the {\it induced} representation of the corresponding
Lie group.

This method is inapplicable both if \fg is infinite-dimensional
and if \fg is a superalgebra, even a finite-dimensional one.

In the infinite-dimensional case the reason is that there is no
correspondence between abstract infinite-dimensional Lie algebras
and infinite-dimensional Lie groups, and no correspondence between
their representations.

For finite dimensional 
Lie superalgebras, the correspondences in question definitely do
exist, but up to now there is no printed text\footnote{%
``Induced representations of a supergroup'' $\cal G$, defined in \cite{S},
are, in fact, representations of the super{\it algebra\/} of $\cal G$,
corresponding to induced representations of $\cal G$.
\DL{sobiralsya zhe upomyanut' Son'ku, a to chto ona -- v alg.
geom. kontexte, tak i chto?!}\VM{Upomjanul tvoju ljubimuju Son'ku. Teper' dovolen?}}
with a {\it definition} of induced
representations of Lie supergroups. Not that it is
difficult to define these creatures,\footnote{%
The topos of sets is not appropriate external universe for the category
of supermanifolds. In particular, coinduced representations of supergroups
can not be defined as sets with some structure (though all modern and future mathematics
is nothing but conservative extension of set theory (any category is a set or class),
not every category is a category of sets with {\it structure}).
There is a universal way to embed any category into appropriate universe
(of contravariant set-valued functors) -- Yoneda's ``point functor''.
For the category of supermanifolds there is an alternative solution: the category
of set-valued functors from the category of finite-dimensional Grassmann superalgebras~\cite{M3}.
The latter approach permits one to extend the category of supermanifolds itself, including
in it infinite-dimensional supermanifolds. Then induced representations of finite-dimensional
supergroups can be defined as {\it internal} objects: Frechet supermanifolds
of ``supersections'' of a vector superbundle  with supersmooth
action of a supergroup on them.}
but the current trend in
``supermath'' is to develop the corresponding issues on ``local''
purely Lie superalgebraic level, ignoring the lifting to the Lie
supergroup level, in cases where such a lifting exists.

In this situation, many people dealing with coinduced
representations of Lie superalgebras, including myself, were
forced to make the necessary calculations introducing 
``odd parameters $\theta$'' and formally imitating the
corresponding calculations for Lie algebra representations.

In~\cite{1} I gave a solid mathematical background to these formal
manipulations for the case most often occurring in practice
(finite-dimensional \fg and $\fg_-$ a Lie subsuperalgebra of \fg),
replacing ``quick and dirty''\footnote{In programmer's slang this means
a quicky bug fix without 
localizing in source code the bug itself, i.e. without real understanding
of the bug's reason.
}
trick with induced
representations by correct manipulations with certain formal power
series.

Here the result of \cite{1} (as well as of I.~Kantor's~\cite{3}
main theorem) is extended for the most general case where both $\fg_-$
and \fh, as well as $V$, may be infinite-dimensional and,
moreover, $\fg_-$ need not be a subsuperalgebra of \fg. I also
give an expression for this action as a ``path integral'' over
paths in a graph generated by the action of $\fg_-$ on \fg. This
expression permits one to write a program, calculating
differential operators representing the 
basis elements of coinduced representations. The program realizing
the corresponding algorithm 
is also written. Since the algorithm for calculating the structure
constants of the simple Lie algebras 
(described in \cite{DG}) is especially simple for simply laced Lie
algebras  $A_n$, $D_n$ and $E_n$ in their Chevalley bases and with
respect to the standard decomposition
$\fg_-=\oplus_{\alpha<0}\fg_\alpha$,
$\fh=\oplus_{\alpha\geq0}\fg_\alpha$, the current (beta) version
of this program calculates the generators and other elements of
the basis of the Lie superalgebra in these cases (in particular,
for \fg-modules dual to Verma modules, the latter being ``free''
highest weight modules). An extension of the program to the
non-simply laced simple Lie algebras, and to simple Lie
superalgebras is being developed.

\section*{Notations and conventions.}
Throughout the paper all Lie superalgebras and vector superspaces
are over a ground field \K of characteristic $0$. 
This restriction is forced mainly by the fact that, in our
explicit constructions, we use formal series in various vector
spaces over \K, whose coefficients contain factorials in
denominators.

Vector superspaces will be considered, if necessary, as {\it
topological} vector superspaces (over \K equipped with the discrete
topology). We tacitly assume that the default topology of a vector
superspace $E$ is discrete again, whereas that of $E^*$ has the
set of all subspaces orthogonal to finite-dimensional subspaces of $E$ as the base of
neighborhoods of zero.

For any Lie superalgebra  \fg, let  $U(\fg)$ be the universal
enveloping superalgebra of \fg. Hereafter if $A$ is an element of  some
superspace (see Ref.~\cite{7} or Appendix~\ref{a1}), we denote the
even (resp. odd) component of $A$ by ${}_\0A$ (resp. ${}_\1A$);
let $[\cdot, \cdot]$ denote  the composition (bracket) in Lie
superalgebras.

The parity of an homogeneous element $X$ of a vector superspace
will be denoted by $|X|$.

We will often define multilinear maps by their values on
homogeneous elements (i.e., either even or odd) without mentioning
this explicitly. The same practice will be applied while checking
identities between multilinear maps, to avoid writing out numerous
``sums over parities'' widening the formulas to an abnormal
extent. The presence of expressions like $(-1)^{|X|\cdot|Y|}$ in a
formula will indicate that the formula is represented in its
``short-hand'' form.

$\mathfrak S_n$ denotes the group of permutations of the set $\{1,\dots,n\}$.

$\hhat E$ denotes the completion of a topological vector
(super)space $E$.

$\overline X$ denotes the closure of a subset $X$ of a topological
space.

\kspan P denotes the \K-linear span of a subset (or of an indexed
family) $P$ of a topological vector (super)space.

If $E$ and $W$ are topological vector superspaces, then $L(E,W)$
denotes the superspace af all {\it continuous} linear maps from $E$ to $W$,
whereas \L E W denotes the superspace $L(E,W)$ equipped with lpco topology,
defined in Sect.~\ref{lpco}.

\bil E,E',W; denotes the superspace of all continuous bilinear maps from
$E\times E'$ to $W$.

Throughout the paper \fg is a Lie superalgebra, \fh a Lie
subsuperalgebra of \fg, $\fg_-$ a subsuperspace of \fg such that
$\fg=\fg_-\oplus\fh$. Depending on the context, $V$ will denote
either an \fh-module, or a superspace.

\section{The induced and coinduced modules over Lie superalgebras.}


Here we reproduce, for the reader's convenience,
general definitions of induced and coinduced modules
as well as some of their properties, see
\cite{9} and \cite{10}.

\subsection{The induced modules.}

Let \fh be a Lie subsuperalgebra of a Lie superalgebra \fg
and 
$V$ an \fh-module.
The \fg-module
\begin{equation}
\label{1.2}
\I V\df U(\fg)\otimes_{U(\fh)} V
\end{equation}
is said to be {\bf induced} by the \fh-module $V$.
As a vector space, \I V is
the quotient space of tensor product $U(\fg)\otimes V$
modulo the relations
\begin{equation}
\label{1.3}
gh\otimes v=g\otimes\rho(h)v\ \ \ (h\in \fh).
\end{equation}

The action of the Lie superalgebra \fg on \I V 
is defined as follows:
\begin{equation}
\label{1.5}
g(u\otimes v)=gu\otimes v\ \ %
(g\in\fg,\ u\otimes v\in \I V),
\end{equation}
where $u\otimes v$ denotes, by abuse of notation, the element 
$u\otimes_{U(\fh)} v$ of \I V.

There is a canonical morphism of \fh-modules
\begin{equation}
\label{1.5'}
\varphi: V\arr \I V|_\fh\ \ \ (\varphi(v)\df1\otimes v)
\end{equation}
such that the pair $(V,\varphi)$ is universal:
for any other pair $(W,\sigma:V\arr W|_\fh)$,
there exists a unique morphism $\theta:\I V\arr W$
of \fg-modules such that $\sigma=\theta\varphi$.

\subsection{The coinduced modules.}
Let \fh, \fg and $V$ be as above.
The \fg-module
\begin{equation}
\label{1.2c}
\T V
\df \Hom_{U(\fh)}(U(\fg),V)
\end{equation}
is said to be {\bf coinduced} (or {\bf produced} in the
original Blattner's terminology \cite{9})
by the \fh-module $V$.
The module \T V as a vector space is
the subspace of $\Hom(U(\fg),V)$
consisting of all \fh-invariant elements $f$: 
\begin{equation}
\label{1.3c}
hf(g)=\sum_{i,j\in\Z_2}(-1)^{ij}{}_if({}_jhg)\ \ \ (h\in \fh).
\end{equation}
(i.e., intertwining maps).

The action
of Lie superalgebra \fg on \T V
is
defined as follows:
\begin{equation}
\label{1.5c}
(gf)(g')=\sum_{i,j,k\in\Z_2}(-1)^{i(j+k)}{}_jf({}_kg'{}_ig)\ \ %
(g,g'\in\fg,\ f\in \T V).
\end{equation}

There is a canonical morphism of \fh-modules
\begin{equation}
\label{1.5c'}
\varphi: \T V|_\fh\arr V\ \ \ (\varphi(f)\df f(1))
\end{equation}
such that the pair $(V,\varphi)$ is universal:
for any other pair $(W,\sigma:W|_\fh\arr V)$,
there exists a unique morphism $\theta:W\arr\T V$
of \fg-modules such that $\sigma=\varphi\theta$.


\subsection{A fine structure of coinduced modules.}

Consider the space \K as the trivial \fh-module.
Given an \fh-module $V$, define the \K-bilinear map
\begin{equation}
\label{module}
\T\K\times\T V\stackrel\cdot\arr\T V
\end{equation}
as follows. Let
\begin{equation}
\label{delta}
\Delta: U(\fg)\arr U(\fg)\otimes U(\fg)
\end{equation}
be a comultiplication in $U(\fg)$, i.e., the only morphism of associative superalgebras that
extends the linear map
\begin{equation}
\label{delta0}
\Delta_0: \fg\arr U(\fg)\otimes U(\fg)\ \ \ (\Delta_0(g)=g\otimes 1+1\otimes g)
\end{equation}
from $\fg\subset U(\fg)$ to the whole $U(\fg)$.

For any $f\in \T\K$ and $s\in \T V$, define $f\cdot s\in\T V$ as follows.
For any $g\in U(g)$, let $\Delta(g)=\sum_ia_i\otimes b_i$. Then:
\begin{equation}
\label{f.s}
(f\cdot s)(g)\df\sum_{i,\eps,\eps'}(-1)^{\eps\eps'}f({}_\eps a_i)\,{}_{\eps'}s(b_i).
\end{equation}

The definition (\ref{f.s}) is just the straightforward extension (via ``sign rule'')
to the case of Lie superalgebras
of the particular case of the corresponding definition from \cite{9}.

The following proposition is the generalization
of Proposition 10 from \cite{9}
to the case of Lie superalgebras.

\begin{Prop}
\label{D.O.}
{\bf a)} The map $(\ref{module})$ 
defines on \T\K
the structure of a supercommutative superalgebra and, for any \fh-module $V$,
turns \T V into an \T\K-module.

\noindent{\bf b)} \fg acts on \T V via differentiations with respect to the just defined structure.
In other words,
\begin{equation}
\label{diff}
g(f\cdot s)=(gf)\cdot s+\sum_{\eps,\eps'\in\Z_2}(-1)^{\eps\eps'}{}_\eps f\cdot({}_{\eps'}gs)
\end{equation}
for any $f\in\T\K$, $\ s\in\T V$, and $\ g\in\fg$.
\end{Prop}

\begin{proof} Heading a) is straightforward. Let us prove b).
Let $g'$ be a homogeneous element in $U(\fg)$ and $\Delta(g')=\sum_ia_i\otimes b_i$.
Then, for any homogeneous $g\in\fg$, we have
\[
\Delta(g'g)=\Delta(g')\Delta(g)=\sum_i((-1)^{|g|\cdot|b_i|}a_i g\otimes b_i+a_i\otimes b_ig)
\]
and
\iffalse
\begin{align*}
&(-1)^{|g|(|f|+|s|+|g'|)}g(f\cdot s)(g')=(f\cdot s)(g'g)\hfill\\
&=\sum_i((-1)^{|g|\cdot|b_i|+|s|(|a_i|+|g|)}f(a_i g)s(b_i)+(-1)^{|a_i|\cdot|s|}f(a_i )s(b_ig))\\
&=\sum_i((-1)^{|g|\cdot|b_i|+|s|(|a_i|+|g|)+|g|(|f|+|a_i|)}gf(a_i )s(b_i)\\
&\ \ \ \ \ \ +(-1)^{|a_i|\cdot|s|+|g|(|s|+|b_i|)}f(a_i )gs(b_i))\\
&=\sum_i((-1)^{|g|\cdot|g'|+|s|(|a_i|+|g|)+|g|\cdot|f|}gf(a_i )s(b_i)\\
&\ \ \ \ \ \ +(-1)^{|a_i|\cdot|s|+|g|(|s|+|b_i|)}f(a_i )gs(b_i))\\
&=\sum_i((-1)^{|g|(|f|+|s|+|g'|)+|s|\cdot|a_i|}gf(a_i )s(b_i)\\
&\ \ \ \ \ \ +(-1)^{|a_i|\cdot|s|+|g|(|s|+|b_i|)}f(a_i )gs(b_i))\\
&=(-1)^{|g|(|f|+|s|+|g'|)}((gf)\cdot s)(g')\\
&\ \ \ \ \ \ +(-1)^{|g|(|f|+|s|+|g'|+|f|\cdot|g|)}(f\cdot(gs))(g')\\
&=(-1)^{|g|(|f|+|s|+|g'|)}(((gf)\cdot s)+(f\cdot(gs)))(g').
\end{align*}
\else
\[
\renewcommand{\arraystretch}{1.4}
\begin{array}{ll}
(-1)^{|g|(|f|+|s|+|g'|)}g(f\cdot s)(g')=(f\cdot s)(g'g)\hfill&\\
=\sum_i((-1)^{|g|\cdot|b_i|+|s|(|a_i|+|g|)}f(a_i g)s(b_i)+(-1)^{|a_i|\cdot|s|}f(a_i )s(b_ig))&\\
=\sum_i((-1)^{|g|\cdot|b_i|+|s|(|a_i|+|g|)+|g|(|f|+|a_i|)}gf(a_i )s(b_i)
+(-1)^{|a_i|\cdot|s|+|g|(|s|+|b_i|)}f(a_i )gs(b_i))&\\
=\sum_i((-1)^{|g|\cdot|g'|+|s|(|a_i|+|g|)+|g|\cdot|f|}gf(a_i )s(b_i)
+(-1)^{|a_i|\cdot|s|+|g|(|s|+|b_i|)}f(a_i )gs(b_i))&\\
=\sum_i((-1)^{|g|(|f|+|s|+|g'|)+|s|\cdot|a_i|}gf(a_i )s(b_i)
+(-1)^{|a_i|\cdot|s|+|g|(|s|+|b_i|)}f(a_i )gs(b_i))&\\
=(-1)^{|g|(|f|+|s|+|g'|)}((gf)\cdot s)(g')
+(-1)^{|g|(|f|+|s|+|g'|+|f|\cdot|g|)}(f\cdot(gs))(g')&\\
=(-1)^{|g|(|f|+|s|+|g'|)}(((gf)\cdot s)+(f\cdot(gs)))(g').&
\end{array}
\]
\fi
The above calculation uses definitions~(\ref{1.5c}), (\ref{f.s}) and
the fact that $|a_i|+|b_i|=|g'|$.
\end{proof}

Proposition~\ref{D.O.} states, roughly, that the action of \fg in
coinduced modules is {\it always} realized by differential
operators of the first order. But it gives no prescription as to
how to find these differential operators. The main purpose of this
paper is just to find {\it explicit} realizations for these
differential operators in a ``natural'' basis in \T V.
The first step is to define a ``natural'' basis.

\subsection{\I V and \T V as superspaces.}
Let $\fg_-$ be a subsuperspace in \fg such that $\fg=\fg_-\oplus\fh$.
We do not require here
that $\fg_-$ be a sub{\it superalgebra} of \fg. Moreover, no condition of finite-dimensionality
is imposed on \fg, \fh or $\fg_-$.

If 
$P=\fam{P_i}{i\in I}$ is a well-ordered basis of $\fg_-$,
then the space \I V is
spanned
by $1\otimes V$ and 
by monomials
\begin{userlabel}
\begin{equation}
\label{1.4}
P_{i_1}^{\alpha_1}\dots P_{i_k}^{\alpha_k}\otimes v\ \ \ (k\in\N\setminus\{0\},\ i_1<\dots<i_k,\ v\in V),
\end{equation}
where 
$\alpha_i=1$ if $P_i$ is an odd 
basis element
and
$\alpha_i$ is a positive integer 
if $P_i$ is an even 
basis element
(The Poincar\'e-Birkhoff-Witt theorem~\cite{11}, \cite{10}).

Hence, every such ordered basis 
defines a linear isomorphism
\begin{equation}
\label{1.4'}
\I V\approx\K[\fam{P_i}{i\in I}]\otimes V\approx S(\fg_-)\otimes V\;,
\end{equation}
\end{userlabel}
where 
$\K[\fam{P_i}{i\in I}]$
is the superspace of supercommuting polynomials,
see Appendix~\ref{a1}, in the $P_i$ (whose number may be infinite).

\medskip
{\bf Remark}. The isomorphism~(\ref{1.4'}) is an
analogue of transition from non-commuting operators
in quantum field theory to operators commuting under
the sign of normal product or $T$-product.

\medskip

The counterpart of (\ref{1.4'}) for coinduced modules is the isomorphism:
\begin{equation}
\label{1.4''}
\T V\approx\Hom(\K[\fam{P_i}{i\in I}],V)\approx\Hom( S(\fg_-) ,V)\;.
\end{equation}

The isomorphisms~(\ref{1.4'}) and (\ref{1.4''}) respect, moreover,
natural \Zplus-filtrations
of \I V and \T V (for details, see \cite{9}).

For our purposes (to find explicit expressions for (co)induced modules)
more appropriate are ``supersymmetrized'' bases in \I V or \T V defined below.

Define the supersymmetric product 
on $\fg_-^k$
with values in $U(\fg)$ as follows:
\begin{equation}
\label{3.3}
g_{1}
\circ
\dots
\circ
g_{k}\df
{1\over{k!}}\sum_{\sigma\in\mathfrak S_k}
(-1)^{\varepsilon(\sigma_{\rm odd})}
g_{\sigma1}\dots g_{\sigma k}\;,
\end{equation}
where $\varepsilon(\sigma_{\rm odd})$
is the parity of permutation of odd elements
among
\cd g{1}{{k}}.

In what follows, we will denote $U(\fg_-,\fg)$ the subsuperspace of $U(\fg)$ spanned by
all supersymmetric products including the empty product $1$. The subsuperspace
$U(\fg_-,\fg)$
clearly is spanned by $1$ and by ``monomials''
\begin{equation}
P_{i_1}\circ\dots\circ P_{i_k}\ \ \ \ (k\in\N\setminus\{0\},\ i_1\le\dots\le i_k).
\end{equation}
If $\fg_-$ is a subsuperalgebra of \fg, then $U(\fg_-,\fg)$ is a subsuperalgebra of $U(\fg)$
coinciding with ``canonical'' subsuperalgebra spanned by non-supersymmetrizied monomials
and isomorphic to $U(\fg_-)$. This justifies the shorthand notation
$U(\fg_-)$ for $U(\fg_-,\fg)$, which will be used in what follows..

The Poincar\'e-Birkhoff-Witt theorem implies that for any direct sum decomposition
$\fg=\fg_-\oplus\fh$ of $\fg$ into pair of subsuperspaces there is an isomorphism
\begin{equation}
\label{PBWiso}
U(\fg_-)\otimes U(\fh)\approx U(\fg)\ \ \ \ (\sum_ig_i\otimes h_i\mapsto\sum g_ih_i).
\end{equation}

The latter isomorphism
implies 
that the superspace \I V is spanned by $1\otimes V$ and by monomials
\begin{equation}
\label{symspan}
P_{i_1}^{\circ\alpha_1}\circ\dots\circ P_{i_k}^{\circ\alpha_k}\otimes v
\ \ \ (k\in\N\setminus\{0\},\ i_1<\dots<i_k,\ v\in V),
\end{equation}
where $g^{\circ k}$ denotes $k$-fold supersymmetric product of $g$.

The isomorphisms~(\ref{1.4'}) and (\ref{1.4''}) arising from the choice of
``supersymmetrized'' basis in $U(\fg)$ also respect, as is easy to see, the
natural \Zplus-filtrations of \I V and \T V.

\section{Linearly topologized superspaces and superalgebras of formal power series.}
Before proceeding further, let us describe some useful topologies
on superspaces arising, in particular, on superalgebras of formal
power series (in possibly infinite number of variables). These
topologies will turn our superspaces into complete topological
superspaces over the field \K with the discrete topology.

\iftrue
The pragmatical reader interested in finite-dimensional case only,
i.e. the case where both the superalgebra \fg and the \fh-module $V$
are finite-dimensional, can skip the most part of this section.
In this case the rudimentary knowledge of elementary properties of
linearly compact spaces (especially of
the superalgebra of formal power series in finite number of indeterminates)
is enough to understand constructions of next sections.

In the general case the superspaces arising are not linearly compact,
so that this section contains all necessary
results from the theory of linearly topologized spaces.
And, for the sake of completeness, a number of more general results
not used in this text, but lacking in canonical sources~\cite{Ko} and~\cite{16},
where this subject is treated.
\else
Though, in what follows, one could do without introducing any topology,
equipping superspaces with suitable topologies will facilitate
the exposition.
One of the benefits of this approach is a free usage of infinite sums over
families of elements of superspaces, interpreted as summable families of elements.
Another benefit is the possibility to extend a (continuous) linear map defined
on a dense subsuperspace of a superspace by continuity to the whole superspace
(if the target superspace is complete).
\fi

\subsection{Linearly topologized superspaces.}
A superspace $E$ equipped with a topology $\tau$ is said to be a
{\bf topological superspace} over \K if it is a topological vector space over \K
and, besides, the decomposition $E={}_\0E\oplus{}_\1E$ is a decomposition
into {\it topological} direct sum (in particular, both ${}_\0E$ and ${}_\1E$
are closed subspaces of $E$).

A Hausdorff topology $\tau$ of a topological superspace $E$ is said to be {\bf linear}
(and $E$ itself as {\bf linearly topologized})
if it has a base of open neighborhoods of~$0$, consisting of subspaces.
Clearly, this is equivalent to the condition that $\tau$ has a basis
of open neighborhoods of~$0$, consisting of sub{\it super\/}spaces.
This is equivalent as well to the condition that both ${}_\0E$ and ${}_\1E$
are linearly topologized vector spaces (the theory of which is
developed in Ch.~2 \ParagraphSign~10 of~\cite{Ko} and,
in more general context of modules over commutative rings,
in Excercises~14-28 in Ch.III, \ParagraphSign~2 of \cite{16}).

\subsection{Discrete superspaces.} A superspace $E$ equipped with the discrete topology is, clearly,
linearly topologized topological superspace briefly called a {\bf discrete superspace}.
In what follows we assume that the Lie algebra \fg, $U(\fg)$, and the \fh-module $V$
are discrete superspaces. The commutative superalgebra
$\S E$ for a discrete superspace $E$ will also be
considered as discrete superspace.
The topological direct sum of any family of discrete superspaces 
is a discrete superspace. Any discrete superspace 
is isomorphic to a topological direct sum
of superspaces $\K^{1|0}$ or $\K^{0|1}$.

\subsection{Cofinite topologies and linearly compact superspaces.}
A linear topology $\tau$ on $V$ is said to be {\bf cofinite}
if it has a base of open neighborhoods of~$0$ consisting of superspaces
of finite codimension. A {\it complete} topological superspace
with cofinite topology is called {\bf linearly compact}.
\begin{Prop}
\label{cofin2lc}
The completion $\wh E$ of a superspace $E$ with cofinite topology
is linearly compact. \eop
\end{Prop}
This justifies the name {\bf linearly precompact} for superspaces
equipped with cofinite topology.

\subsubsection*{Example 1: superspaces $E^*$.} If $E$ is a discrete topological superspace,
then the superspace $E^*$ equipped with weak topology,
i.e., arising from duality between $E^*$ and $E$, is linearly compact.
Recall that a base of neighborhoods of~$0$ of~$E^*$ (both open and closed) in weak topology
is formed by subsuperspaces $K^\perp\subset E^*$ of linear forms annihilating some finite-dimensional
subsuperspace $K$ of $E$.

In what follows superspaces dual to discrete ones will be always considered as
linearly compact topological superspaces. Any topological
product of linearly compact superspaces is linearly compact and
any linearly compact superspace is isomorphic to a topological product
of superspaces $\K^{1|0}$ or $\K^{0|1}$. In particular, any linearly compact
superspace is isomorphic (as topological superspace) to a superspace $E^*$ for some
discrete superspace $E$.

\subsubsection*{Example 2: the cofinite topology generated by a basis.}
\label{tauP}
Let $P=\fami{P_i}$ be a basis of a superspace $E$. Define
a cofinite topology $\tau P$ on $E$ having a base of open neighborhoods of~$0$
consisting of all subsuperspaces
\begin{equation}
\label{basis2noods}
E_F\df\kspan{\fam{P_i}{i\in I\setminus F}},
\end{equation}
where $F$ is a finite subset of $I$.
For infinite-dimensional $E$ the topology $\tau P$ is strictly {\it weaker}
than the weak topology of $E$ (i.e., coming from duality between $E$ and $E^*$).
Recall~(\cite{Ko}) that a basis of open neighborhoods of~$0$ in weak topology
is formed by {\it all} subsuperspaces of $E$ of finite codimension.
For example, the hyperplane $\la\fam{P_i-P_j}{i,j\in I}$ is open in the weak
topology of $E$ but is not open in the topology $\tau P$.

{\small
Let $\wh E_w$, resp. $\wh E_{\tau P}$ be the completion of $E$ with respect to the weak topology,
resp. with respect to the topology $\tau P$.
Both $\wh E_w$ and $\wh E_{\tau P}$ are linearly compact due to Prop.~\ref{cofin2lc} and
there is continuous linear map $\wh\Id:\wh E_w\arr\wh E_{\tau P}$.
We will see below that $\wh E_{\tau P}$ is isomorphic (non-canonically) to $E^*$,
whereas $\wh E_w$ is known~(see~\cite{Ko}) to be isomorphic to $(E^*)_{\text{disc}}^*$,
where $(E^*)_{\text{disc}}$ is the superspace $E^*$ with discrete topology.
If $E$ is infinite-dimensional, then
$\dim (E^*)_{\text{disc}}^*=2^{\dim E^*}$, so that
the superspace $(E^*)_{\text{disc}}^*$ is ``bigger'' than $E^*$ and the map $\wh\Id$
is to have non-zero kernel of dimension $(E^*)_{\text{disc}}^*$.}

\begin{Prop}
\noindent{\bf(a)} The product of any family of linearly precompact superspaces is linearly precompact;

\noindent{\bf(b)} The Image $f(C)$ of any linearly precompact subsuperspace $C$ of a linearly
topologized superspace $E$ under any continuous linear map $f$ is linearly precompact;

\noindent{\bf(c)} Any linearly precompact subsuperspace $C$ of a topological direct sum
$\oplus_{i\in I}E_i$
of a family of linearly topologized superspaces is contained in the direct sum of some finite subfamily.
\end{Prop}\eop

\subsection{Topological superspaces \L E W and $E\botimes W$.}
\label{lpco}
Other types of topological superspaces, which will occur,
are $L(E,W)$ and $E\otimes W$ 
for discrete or linearly compact superspaces $E$ and $W$, equipped with appropriate
topologies.
But we will give the definitions of these topologies
for the general case of arbitrary linear topologies on $E$ and $W$.
After all, the class of linearly topologized superspaces consisting
of all discrete and all linearly compact superspaces only, is not closed
with respect to such categorical operations as $\oplus$, $\otimes$, $\L \cdot\cdot$, etc.

\iffalse
Note
that in the category of all linearly topologized
superspaces, similarly to the category of locally convex vector (super)spaces,
there is no topologies on \L E W and $E\otimes W$, which could pretend to be
``the'' natural ones in all circumstances. There are several definitions instead,
each of which has its own region of applicability: depending of the context
one of these definition may be more ``natural'' than the others.
In what follows we will restrict ourselves with only {\it one}
topology on \L E W and {\it two} topologies on $E\otimes W$.
\else
From~\cite{Ko} one can conclude that the category \LVect of linearly topologized spaces
is both complete and cocomplete, because it has arbitrary direct sums, products,
kernels and cokernels. What about categorical tensor product and internal $\Hom$ functor?
Recall that a bifunctor $(E,W)\mapsto E\botimes W$, resp. $(E,W)\mapsto\mathcal Hom(E,W)$ is called
a {\bf tensor product} functor, resp.~{\bf internal Hom} functor,
if there exists a natural isomorphism
\begin{equation}
\label{tensor}
\bil E,W,V;\approx L(E\botimes W, V),
\end{equation}
resp.
\begin{equation}
\bil E,W,V;\approx L(E,\mathcal Hom(W, V)),
\end{equation}
where \bil E,W,V; is the set of continuous bilinear maps from $E\times W$ to $V$.
We will see that categorical tensor product functor $\botimes$ exists in \LVect.
As to $\mathcal Hom$, it seems to be non-existent on the whole category \LVect.
But there are 2 full coreflective subcategories of \LVect, where internal $\Hom$ does exist.
\fi

\subsubsection*{Linearly precompact-open topology on \L E W.}
Define a linear topology on the superspace $L(E,W)$ 
of continuous linear maps from $E$ to $W$ as follows.
For $C\subset E$ and $O\subset W$, let $L(C,O)$ be the subset of $L(E,W)$, 
consisting of all $f$ such that $f(C)\subset O$.
The superspace $L(E,W)$ 
will be equipped with a linear topology
having as a base of open neighborhoods of~$0$ all subsuperspaces $L(C,O)$,
where $C$ is a 
linearly precompact (closed) subsuperspace of $E$ and $O$ is an open
subsuperspace in $W$. This topology will be called {\bf linearly
precompact-open} or, briefly, {\bf lpco} topology on $L(E,W)$. 
In what follows we will use lpco topology as default one and
the superspace $L(E,W)$ equipped with this topology will be denoted \L E W.

\subsubsection*{Remark.}
Replacing in the above definition linearly precompact subsuperspaces of $E$
by linearly compact ones one gets a {\bf linearly compact-open} (or, briefly, {\bf lco}) topology
on $L(E,W)$. 
This topology is, generally speaking, {\it weaker} than lpco topology,
but if $E$ is complete, then any closed linearly precompact subsuperspace $C\subset E$
is 
linearly compact, so that in this case both topologies coincide.
The superspace $L(E,W)$ equipped with lco topology will be denoted $\L E W_{lco}$.
\\

If $W$ is discrete, then, clearly, the family of subsuperspaces $L(C,0)$ form a base
of open neighborhoods of~$0$ in \L E W. In particular, if $W=\K$ and $E$ is discrete,
then  the lco topology on $E^*$ coincides with the weak topology of $E^*$.

It is clear as well that
if $E=V^*$, then the family of subsuperspaces $L(E,O)$ form a base
of open neighborhoods of~$0$ in \L E W. In particular, if $W=\K$,
then the lco topology on $V^{**}\approx V$ coincides with the discrete topology of $V$.

\begin{Prop}
Let $E$ and $W$ be discrete superspaces. The map
\begin{equation}
*:\L E W\arr\L{W^*}{E^*}\ \ \ \ (f\mapsto f^*)
\end{equation}
is a linear homeomorphism.
\end{Prop}

\begin{proof}
The map $*$ is known to be a linear isomorphism (see~\cite{Ko}),
so it only needs to prove that lpco topologies on \L E W and
\L{W^*}{E^*} are isomorphic. Linearly compact subsuperspaces of a
discrete superspace $E$ are exactly finite-dimensional
subsuperspaces $K$ (see~\cite{Ko}). Given such a superspace
$K\subset E$, we have: 
\[
\renewcommand{\arraystretch}{1.4}
\begin{array}{lll}
*(L(K,0))&=&\{\,f^*\mid f(K)=0\,\}\\
&=&\{\,f^*\mid\varphi(f(K))=0\text{ for all }\varphi\in W^*\,\}\\
&=&\{\,f^*\mid(f^*\varphi)(K)=0\text{ for all }\varphi\in W^*\,\}\\
&=&L(W^*,K^\perp)
\end{array}
\]
And, by the definition of weak topology, the family of spaces $K^\perp$
is a base of open neighborhoods of~$0$ in $E^*$, hence
the family of superspaces $L(W^*,K^\perp)$ is a base of neighborhoods
in \L {W^*}{E^*} in lpco topology.
\end{proof}

The completeness of the superspace \L E W for discrete $E$ and $W$
follows from general Prop.~\ref{compl} below. To formulate this proposition
we need the notion of {\it linerly precompactly generated superspaces\/}.

\subsubsection*{Linerly (pre)compactly generated superspaces.}
Call a linearly topologized superspace $E$ {\bf linearly (pre)\-com\-pact\-ly generated}
or, briefly, {\bf l(p)cg} if any subsuperspace $U$ of $E$ is open if its intersection
with any linearly (pre)compact subsuperspace $C$ is open in $C$.

Any discrete or linearly compact superspace is, clearly, both lcg and lpcg.
Any linearly precompact superspace is lpcg; it is lcg iff it is complete.

The main property of l(p)cg superspaces is:
\begin{Prop}
\label{lpcgmain}
A linear map $f:E\arr W$ of a l(p)cg $E$ into a linearly topologized superspace $W$
is continuous iff for any linearly (pre)compact subsuperspace $C$
of $E$ the restriction $f|_C$ is continuous.
\end{Prop}\eop

\begin{Prop}
\label{compl}
If $E$ is lpcg 
{\rm(resp.} lcg\/{\rm)}
and $W$ is complete, then \L E W
{\rm(resp.} $\L E W_{lco}${\rm)}
is complete.
\end{Prop}
\begin{proof}
The statement follows immediately from~\ref{lpcgmain} and
Theorem~2 in Chapter~10 \ParagraphSign 3 of~\cite{B1}.\VM{Nomer Teoremy i pargrafa -- fal'shivye, najti pravil'nye!!!}
\end{proof}


\small
Denote \LVect the category of linearly topologized vector superspaces,
\PCGLVect (resp. \CGLVect) the full subcategory of \LVect, consisting
of all lpcg superspaces (resp. of all lcg superspaces).

\begin{Prop}
\label{adjoint}
The inclusion functor $\PCGLVect\subset\LVect$ {\rm(resp.} $\CGLVect\subset\LVect$\/{\rm)}
has right adjoint
\begin{equation}
K:\LVect\arr\PCGLVect\ \ \ \ (\text{\rm resp.} K:\LVect\arr\CGLVect).
\end{equation}
\end{Prop}
\begin{Cor}
The categories \PCGLVect and \CGLVect are complete and cocomplete.
\end{Cor}

\begin{Th}
The category \PCGLVect  {\rm(resp. \CGLVect)} have internal $\Hom$ functor.
This functor is defined as $\mathcal Hom(E,W):= K\L E W$ 
{\rm(resp. as $\mathcal Hom(E,W)\df K\L E W_{lco}$).} 
\end{Th}

\begin{Prop}
\begin{equation}
(\bigoplus_{i\in I}KE_i)^*\approx\prod_{i\in I}(KE_i)^*;
\ \ \ \ (\prod_{i\in I}KE_i)^*\approx\bigoplus_{i\in I}(KE_i)^*;
\end{equation}
\end{Prop}

Considering, as usual, a linearly topologized superspace $E$
as a subsuperspace of its double topological dual $E^{**}$,
we call $E$ {\bf reflexive} 
if $E^{**}\approx E$ (as topological superspaces).

We have seen that both discrete superspaces and dual to them
(i.e., linearly compact) superspaces are reflexive.
The next Corollary states that the class of reflexive l(p)cg is big enough:

\begin{Cor}
\label{reflsum}
Any topological direct sum or product of reflexive l(p)cg superspaces is reflexive.
\end{Cor}

Props.~\ref{adjoint}--\ref{reflsum} above
(which proofs and more details will be published soon in~\cite{M2}) show
that both categories \PCGLVect and \CGLVect are counterparts
of the category of Kelley vector spaces (\cite {FrJa},\cite{Po}).

In fact, the direct counterpart is the subcategory \CGLVect of linearly compactly generated
superspaces, but 
bigger subcategory \PCGLVect 
is more flexible in some respects.
\normalsize

\subsubsection*{Topologies on $E\otimes W$.}
We will be interested here only in topologies on superspaces of the type
$E^*\otimes W$ with both $E$ and $W$ being discrete superspaces.
Our aim is to ``approximate'' somehow the superspace \L E W by the superspace $E^*\otimes W$.
But definitions and theorems of this section are given in their natural generality --- for
arbitrary linearly topologized spaces.

From purely categorical point of view ``the natural'' topology on the superspace
$E\otimes W$ should be chosen in such a way,
that it permits one to represent {\sl bi\/}linear continuous maps from $E\times W$
(which are {\sl not} morphisms of the category of topological superspaces)
by continuous {\it linear} maps from $E\otimes W$. In other words
it permits one to internalize definitions of \K-algebras, $R$-modules, etc.
entirely inside the category of linearly topologized superspaces, not using
external maps such as bilinear or, more generally, multilinear ones.
So we will define the {\bf default}\footnote{
I would prefer the name ``tensor product topology'' but this name is already
reserved by Bourbaki  :(( for another linear topology on $E\otimes W$
(the definition in Exc.~28 in Ch.III, \ParagraphSign~2 of \cite{16}, reproduced below).

By the way, the counterpart of default topology for the tensor product of locally convex spaces
is called, strangely enough, {\it projective} topology (see e.g. \cite{Sh}), being in fact
a kind of {\it inductive} one.}
topology on $E\otimes W$
as the strongest linear topology such that the canonical bilinear map
\[
E\times W\arr E\otimes W\ \ \ \ (\,(e,w)\mapsto e\otimes w)
\]
is continuous.
If such linear topology exists, then $E\otimes W$ equipped with this topology is clearly
the categorical tensor product $E\botimes W$, i.e., there is a natural isomorphism~(\ref{tensor}).
\begin{Prop}
For any linearly topologized superspaces $E$ and $W$, the default topology on $E\otimes W$ exists.
A subsuperspace $L\subset E\otimes W$ is 
open
in the default topology iff
for any $e\in E$ and $w\in W$ there exist 
an open subsuperspace $U$ of $E$
and 
an open subsuperspace $U'$ of $W$
such that
\[
e\otimes U'+U\otimes w+U\otimes U'\subset L.\ \eop
\]
\end{Prop}

\iftrue 
In what follows $E\botimes W$ denotes the tensor product of $E$ and $W$
equipped with default topology, whereas $E\hbotimes W$ denotes the completion
of $E\botimes W$ with respect to the default topology.
\fi

Another good property of categorical tensor product is that it respects
arbitrary direct sums:
\begin{Prop}
\label{tensor-directsum}
For any family \fam{W_i}{i\in I} of linearly topologized superspaces
and any linearly topologized superspace $E$ there is a natural isomorphism
\begin{equation}
E\botimes\bigoplus_{i\in I}W_i\approx\bigoplus_{i\in I}\left(E\botimes W_i\right).
\end{equation}
\end{Prop}
But just this ``good'' property prevents us from using categorical tensor product
to approximate the superspace of linear maps. Because Prop.~\ref{tensor-directsum}
implies immediately that, for discrete superspaces $E$ and $W$,
the superspace $E^*\botimes W$ is complete (as direct sum of a number of copies of
the complete superspace $E^*$). And the image of this superspace in \L E W
(see Prop.~\ref{WhotimesE} below)
consists of the maps of finite rank only, so there is no way to approximate maps
of infinite rank via completing $E^*\botimes W$ -- it is complete already.

So we are forced to consider
another
topology on $E\otimes W$ (defined in Exc.~28 in Ch.III, \ParagraphSign~2 of \cite{16}).
It has a base of open neighborhoods of~$0$ consisting of all superspaces
$U\otimes W+E\otimes U'$, where $U$ is an open subsuperspace of $E$
and $U'$ is an open subsuperspace of $W$. Bourbaki calls this topology
the {\bf tensor product topology} but, in my opinion, the name
{\bf projective} topology is more appropriate, because the completion
of $E\otimes W$ in this topology coincides with the projective limit
$\lim (E/U_\alpha\otimes W/U'_\beta)$, where $U_\alpha$, resp. $U'_\beta$ run over
open subspaces of $E$, resp. of $W$.
We will denote $E\otimesp W$ or, by abuse of notations, $E\otimes W$ the tensor product of $E$ and $W$
equipped with the projective topology.

We will need this definition only when $W$ is discrete and $E$ is either discrete or linearly compact
(or vice versa).
In the first case, $E\otimesp W$ is discrete again, hence, complete.
In the second case ($E=V^*$), subsuperspaces $U\otimes W$ form a base of open neighborhoods of~$0$
and
$E\otimesp W$ is not, generally speaking, either discrete,
or linearly compact, or complete.

The completion of $E\otimesp W$ will be denoted $E\hotimesp W$.
The importance of this topology is caused by the following

\begin{Prop}
\label{WhotimesE}
Let $E$ and $W$ be 
linearly topologized
superspaces. Let $\alpha$ be the canonical linear inclusion map
\begin{equation}
\label{otimes2hom}
\alpha:W\otimes E^*\arr\L E W
\ \ \ \ %
(\alpha(w\otimes f)e\df wf(e)\text{ for any }e\in E,\ w\in W, \ f\in E^*).
\end{equation}
Then:

\noindent{\bf(a)} The set of continuous linear maps from \L E W belonging to image of $\alpha$
coincides with the set of continuous linear maps of {\bf finite rank},
i.e. maps with finite-dimensional image.

\noindent{\bf(b)} $\alpha$ is continuous in projective topology on $W\otimes E^*$\/{\rm;}

\noindent{\bf(c)} Moreover, the projective topology on $W\otimes E^*$ coincides with the topology
induced from the lpco topology of \L E W along $\alpha$.

\noindent{\bf(d)} 
The set $\alpha(W\otimes E^*)$ is dense in \L E W.
If, moreover, $E$ is lpcg and $W$ is complete
then $\alpha$ extends by continuity to an isomorphism:
\begin{equation}
\label{otimes2hom^}
\hhat\alpha: W\hotimes E^*\approx\L E W.
\end{equation}
\end{Prop}
In what follows, we will identify topological superspace \L E W with $W\hotimes E^*$ (or with $E^*\hotimes W$)
via the isomorphism~(\ref{otimes2hom^}), which clearly exists if both $E$ and $W$ are discrete
superspaces.

\begin{proof}
%
Heading (a) is evident.

The inverse image under $\alpha$ of a neighborhood $L(C,U)$ of~$0$ in \L E W
is, clearly, $W\otimes C^\perp+U\otimes E^*$, so that $\alpha$ is continuous and heading (b) is valid.
But the set of these inverse images
is a base of neighborhoods of~$0$ in $W\otimes E^*$ by definition of the tensor product topology,
so (c) is also valid.

Prove now heading (d).
Let $f\in\L E W$. We are to prove that for any linearly precompact closed subsuperspace $C$ of $E$
and any open subsuperspace $U$ of $W$ there exists $f_0\in\L E W$, such that
$f-f_0\in L(C,U)$ and $f_0$ is of finite rank.

Consider the subsuperspace $C\bigcap f^{-1}(U)\subset E$.
It is of finite codimension in $C$,
because $f^{-1}(U)$ is open and $C$ is linearly precompact.
Let $e_1,\dots,e_n$ be a basis in an algebraic complement of $C\bigcap f^{-1}(U)$ in $C$.
The \K-span \kspan{e_1,\dots,e_n} is closed in $C$ and $E$ and its intersection with $f^{-1}(U)$ is $0$.
The counterpart of Hahn--Banach theorem for linearly topologized superspaces
(see ``Erweiterungszatz'' on p.~89 of \cite{Ko}) implies, that there exist continuous linear functionals
$\varphi_1,\dots\varphi_n\in E^*$ such that for any $1\le i,j\le n$ one has:
$\varphi_i(e_j)=\delta_{ij}$ and, moreover, $f^{-1}(U)\subset \text{Ker}\,\varphi_i$ for any $i$.
One easily checks then that the finite rank map
\[
f_0:=\alpha(\sum_if(e_i)\otimes\varphi_i)
\]
satisfies the desired condition.

The second part of heading~(d) follows directly from Prop.~\ref{compl}.
\end{proof}

\subsection*{Remark.}
The theory of linearly topologized superspaces is a ``toy'' counterpart of
the theory of locally convex spaces (at least of this part of the latter theory,
which is based on Hahn-Banach theorem). But the first part of heading~(d) of Prop.~\ref{WhotimesE}
is the only result I know of, which has a ``toy'' proof, whereas its locally convex counterpart
(denseness of $E^*\otimes V$ in \L E V equipped with compact-open topology)
is a very difficult problem to solve. In fact, it was yet the open problem in 1971
(see e.g.~\cite{Sh}). I do not know whether this conjecture was proved since then.

\subsection{Topological bases.}
The next thing to do is, given some basis in a (discrete)
superspace $E$, to describe some kind of ``dual basis'' in $E^*$.
Unfortunately, if $E$ is infinite-dimensional, there is no {\it
algebraic} basis which may be called dual to the original basis in
$E$. Just because the algebraic dimension of $E^*$ is then
strictly {\it bigger} than that of $E$:
$\dim E^*=
2^{\dim E}$.
Nevertheless, in this situation
there is an adequate replacement of algebraic basis, namely the
{\it topological} basis of $E^*$.

Let $E$ be a topological superspace. A family $P=\fam{P_i}{i\in I}$
of elements of $E$ is called a {\bf topological basis} of $E$
(or {\bf continuous basis} of $E$ in the terminology of~\cite{Ko})
if any $e\in E$ has the unique representation $e=\sum_{i\in I} c_iP_i$ as a {\it topological} sum
and, besides, for any $i\in E$, the map $\pi_i:E\arr\K\ (e\mapsto c_i)$
is continuous. 

For a discrete superspace $E$, the notions of
topological and algebraic bases clearly coincide.

On the contrary, in $E^*$, there does not exists an algebraic basis
which is, simultaneously, a topological one. Nevertheless (\cite{Ko}):
\begin{Prop}
\label{dualbasis}
Let $P=\fam{P_i}{i\in I}$ be a basis in $E$.

Let the family
$X=\fam{X^i}{i\in I}$
in $E^*$ be defined by the equations: 
\begin{equation}
\label{dualbasis'}
X^i(P_j)=\delta^i{}_j\;.
\end{equation}
Then $X$ is a topological basis in $E^*$. \eop
\end{Prop}
The topological basis $X$ will be called {\bf dual} to the basis $P$;
in what follows we will write simply ``dual basis'' instead of
``dual topological basis''.

S.~K\"othe claimed in~\cite{Ko} that it is not known, for a general 
linearly topologized space $E$, whether or not there exists a topological basis of $E$.
We will only use topological bases for the superspaces dual to discrete superspaces,
and, more generally, for the superspaces $\L E W\approx W\hotimes E^*$ for discrete $E$ and $W$.

For $E^*$ such bases always exist due to Prop.~\ref{dualbasis}, and, moreover,~(\cite{Ko}):
\begin{Prop}
Any topological basis of $E^*$ is dual to some basis of $E$. \eop
\end{Prop}
And what about \L E W and $W\hotimes E^*$? The following Proposition generalizes Prop.~\ref{dualbasis}
to this case.
\begin{Prop}
Let $\fam{P_i}{i\in I}$ be a basis in $E$, $\fam{X^i}{i\in I}$ be a dual basis in $E^*$
and let $\fam{P'_k}{k\in K}$ be a basis in $W$.

\noindent{\bf(a)} Let the family $X=\fam{X^i_k:E\arr W}{i\in I,k\in K}$ be defined by equations:
\begin{equation}
\label{hombasis'}
X^i_k(P_j)=\delta^i{}_jP'_k\;.
\end{equation}
Then $X$ is a topological basis in \L E W. The inverse image
of this basis in $W\hotimes E^*$ under the isomorphism $\hhat\alpha:W\hotimes E^*\approx\L E W$
{\rm(see~(\ref{otimes2hom^}) in Prop.~\ref{WhotimesE}~(d))} is a topological basis
\begin{equation}
\wh\alpha^{-1}(X)=\fam{\wh\alpha^{-1}(X^i_k)=P'_k\otimes X^i}{i\in I,k\in K}.
\end{equation}

\noindent{\bf(b)} The formal sum
\[
\sum_{i,k}c_i^kX^i_k\ \ \ \ (\text{resp. }\sum_{i,k}c_i^kP'_k\otimes X^i)
\]
is a topological sum representing some element of \L E W {\rm(resp. of $W\hotimes E^*$)}
iff, for any $i\in I$, there exists only finite number of $k\in K$ such that $c_i^k\neq0$. \eop
\end{Prop}

\begin{Cor}
\label{virtbasis}
\noindent{\bf(a)}
For any element $x\in W\hotimes E^*$ there exists a unique family \fami{w_i}
of elements of $W$
{\rm(}resp. a unique family \fami{f^i} of elements of $E^*${\/\rm)}
such that $x$ represents as a topological sum
\begin{equation}
x=\sum_{i\in I}w_i\otimes X^i=\sum_{i\in I}P'_i\otimes f^i.
\end{equation}
\noindent{\bf(b)} For any family \fami{w_i} of elements of $W$ the topological sum
$\sum_{i\in I}w_i\otimes X^i\in W\hotimes E^*$ exists. \eop
\end{Cor}

\subsection*{Remark.} Heading (b) of Cor.~\ref{virtbasis} permits one
to consider informally a basis  $\fam{X^i}{i\in I}$ 
in $E^*$ as a ``virtual'' basis of $W\hotimes E^*\approx\L E W$
``with coefficients in $W$''.
Meditate a bit over this informal interpretation :).

\subsection{The superalgebra $S(E)^*$.}
\label{SE*}
Consider a superspace $E$ as a commutative Lie superalgebra (with the zero bracket),
let $S(E)^*$ be an $E$-module coinduced from trivial 1-dimensional representation
of the trivial subalgebra $\fh=0$ of $E$. Then $S(E)^*$
has a natural structure of a commutative associative superalgebra with unit element
and multiplication defined by eqs.~(\ref{f.s}). How this superalgebra looks like?

First of all, the topological superspace $E^*$ may be identified with a closed subsuperspace
of $S(E^*)$ via continuous injective map
\begin{equation}
\label{Estar2SEstar}
\iota:E^*=S^1(E)^*\harr\bigoplus_{n\in\N}S^n(E)^*\arr\prod_{n\in\N}S^n(E)^*
\approx(\bigoplus_{n\in\N}S^n(E))^*=S(E)^*.
\end{equation}
In particular, $E^*$, considered as a subsuperspace of $S(E)^*$, coincides with
the annihilator of the subsuperspace $\mathop{\bigoplus}\limits_{n\neq 1}S^n(E)$ of $S(E)$.

In what follows we identify $E^*$ with a subsuperspace of $S(E^*)$,
writing just $f$ instead of $\iota(f)$ for $f\in E^*$.

\begin{Prop}
\label{p6}
For any homogeneous $f_1,\dots,f_n\in E^*\subset S(E)^*$ and any homogeneous $e_1,\dots,e_m\in E$,
we have
\begin{equation}
\label{3.8}
(f_1\dots f_m)(e_n\dots e_1)=
\delta_{mn}\sum_{\sigma\in\mathfrak S_n}
(-1)^{\varepsilon(\sigma_{\rm odd})}
f_{\sigma(n)}(e_n)\dots
f_{\sigma(1)}(e_1).
\end{equation}

\end{Prop}

\begin{proof}
For a finite well ordered set $S$ and a natural number $i$ such that $i<\text{card}(S)$,
the image of $i$ under the only order preserving isomorphism
$\text{card}(S)\arr S$ will be denoted $S_i$.
For any nontrivial partition of a well ordered set $n=\{\,0,\dots,n-1\,\}=S\cup S'$,
where $S'=n\setminus S$,
denote by $\sigma(S)$ the permutation
\[
(0,\dots,n-1)\mapsto(S_0,\dots,S_{\text{card}(S)-1},S'_0,\dots,S'_{\text{card}(S')-1}).
\]
For trivial partitions ($S=\emptyset$ or $S=n$), let $\sigma(S)$ be the identity permutation $\One_n$.

In these notations, we can express the action of comultiplication $\Delta$
on any element $e_0\dots e_{n-1}$ of $S^n(E)$ as follows: 
\begin{equation}
\label{Delta}
\Delta(e_0\dots e_{n-1})=\sum_{S\subset n}(-1)^{\eps\sigma(S)_{\text{odd}}}
e_{S_0}\dots e_{S_{\text{card}(S)-1}}\otimes e_{S'_0}\dots e_{S'_{\text{card}(S')-1}}.
\end{equation}
The expression~(\ref{Delta}) clearly follows from the fact that $\Delta$ is morphism of superalgebras:
\[
\Delta(e_0\dots e_{n-1})=\Delta(e_0)\dots\Delta(e_{n-1})=
(e_0\otimes 1+1\otimes e_0)\dots(e_{n-1}\otimes 1+1\otimes e_{n-1}).
\]
Now eq.~(\ref{3.8}) can be proved by induction on $n$
if one applies the definition~(\ref{f.s}) to $f_1\dots f_n=f_1(f_2\dots f_n)$.
\end{proof}

\begin{Cor}
\label{SE*=K[X]}
The choice of a basis $P$ in $E$ establishes an isomorphism of topological superspaces
\begin{equation}
\label{3.9}
S(E)^*\approx\K[[X]].
\end{equation}
Moreover, this isomorphism respects multiplication, i.e. is an isomorphism of superalgebras.
\end{Cor}
We will identify $S(E)^*$ and $\K[[X]]$ via this isomorphism.
In fact the isomorphism~(\ref{3.9}) implies that
$S(E)^*$ is a complete {\it topological} superalgebra,
because $\K[[X]]$ is (see Exc.~(22) to \ParagraphSign~2 of Ch.~3 of~\cite{16}).

The following proposition allows one to explicitly calculate the expression for
\T\M 
for any $\M\in \fg$ in any given pair of homogeneous bases:
\fam{P_i}{i\in I} 
in $\fg_-$
and the dual basis
\fam{X^i}{i\in I} 
in $\fg_-^*$:
\begin{Prop}
\label{p8}
The isomorphism~$(\ref{3.9})$ establishes the
``Fourier transform'' of the operators $\partial\over\partial X^i$
and left multiplications by $P_i$ in $S(E)$:
\begin{equation}
\label{3.18}
P_i{}^*=(-1)^{|X^i|}{\partial\over\partial X^i}\;,\ \
\left({\partial\over\partial P_i}\right)^*=X^i,
\end{equation}
where we have identified $P_i{}^*$ and
$\left({\partial\over\partial P_i}\right)^*$ with their images
under the isomorphism~$(\ref{3.9})$. 
\end{Prop}

\begin{proof}
The validity of relations~(\ref{3.18}) on 
the dense subsuperspace
$\K[X]\approx S(\kspan{X})\subset S(E)^*$
is verified straightforwardly
by applying Eq.~(\ref{3.8}).
Clearly, the operators
$P_i{}^*$ and
$\left({\partial\over\partial P_i}\right)^*$
are continuous, so that
eqs.~(\ref{3.18}) can be extended to hold on
the whole $S(E)^*$. 
\end{proof}

Define now the bilinear map
\begin{equation}
\label{3.12}
\sar\Sstar E\hotimes S(E),\cD,\diff{\Sstar E};\;,
\end{equation}
where \diff{\Sstar E} is the space of differential operators
of finite order in $\Sstar E$.
This map participates in the expression for  the action of \fg on \T V.
The map \cD is defined as follows:
\begin{equation}
\label{3.13}
\cD(f\otimes e)\varphi\df f\cdot(e*\varphi)\ \ %
(\varphi,f\in\Sstar E,\ \ e\in E).
\end{equation}
Here $e^*$ is the operator dual to the operator of
left multiplication by $e$ in $S(E)$. Proposition~\ref{p8} above 
implies that if $e=e_1\dots e_n$, where $e_i\in E$,
then
\[
e^*=(e_1\dots e_n)^*=e_1^*\dots e_n^*
\]
is a product of differentiations in $\Sstar E$.

In particular, $P_i^*=(-1)^{|X^i|}{\partial\over{\partial X^i}}$, so that
one has the following expression for the action of \cD in the basis
$P=\fam{P_i}{i\in I}$ on $E$ and the dual basis $X=\fam{X^i}{i\in I}$ on $E^*$:
\begin{equation}
\label{3.13expl}
\cD(f(X)\otimes P_{i_1}\dots P_{i_n})=(-1)^{|X^{i_1}|+\dots+|X^{i_n}|}
f(X){\partial\over\partial X^{i_1}}\dots{\partial\over\partial X^{i_n}} 
\end{equation}
or, symbolically:
\begin{equation}
\label{3.13explsym}
\cD(\varphi(X,P))=\varphi(X,(-1)^{|X|}{\partial\over\partial X}).
\end{equation}

Another  bilinear map
\begin{equation}
\label{3.12 infty}
\sar S(E)\hotimes\Sstar E,\cD_\infty,\diffinf{S(E)};\;,
\end{equation}
where \diffinf{S(E)} is the space of differential operators
of {\it infinite\/} order with {\it polinomial\/} coefficients on $\hhat S(E)$,
participates in the expression for the action of \fg on \I V.
The map $\cD_\infty$ is defined as follows:
\begin{equation}
\label{3.13 infty}
\cD_\infty(e\otimes f)v\df e\cdot (f^*v)\ \
(e,v\in S(E),\ \ f\in \Sstar E).
\end{equation}

The counterparts of (\ref{3.13expl}) and (\ref{3.13explsym}) are:
\begin{equation}
\label{3.13explinf}
\cD_\infty(P_{i_1}\dots P_{i_n}\otimes f(X))=P_{i_1}\dots P_{i_n} f({\partial\over\partial P}) 
\end{equation}
or, symbolically:
\begin{equation}
\label{3.13explsyminf}
\cD_\infty(\varphi(P,X))=\varphi(P,{\partial\over\partial P}).
\end{equation}

The reader who does not believe in existence of differential operators of infinite order
may interpret~(\ref{3.13 infty}) just as {\it definition} of the space \diffinf{S(E)} :-).

As to the space \diff{\Sstar E} (the existence of which is of no doubt), one can prove
that \cD is a \K-linear {\it isomorphism\/}.
We need, in fact, only the following easy result:
\begin{Prop}
\label{diffs are good}
The restriction of \cD on the subsuperspace $\Sstar E\hotimes E$
establishes an isomorphism of this subsuperspace with the superspace
of all differentiations of the superalgebra $\Sstar E\approx S(E)^*$.
\end{Prop}

\begin{proof}
For any differentiation $d$, we have $d=\sum_id(X^i){\partial\over\partial X^i}$.
\end{proof}

\section{Main theorem.}
Let
\begin{equation}
\label{a3.1}
\fg=\fg_-\oplus\fh
\end{equation}
be a decomposition of a Lie superalgebra \fg into the
direct sum of subsuperspaces.

In what follows we will identify via Props.~\ref{WhotimesE} (d) 
the topological superspace
\begin{equation}
\label{a3.13}
\T V\approx\L{\S{\fg_-}}V
\end{equation}
with the topological superspace
\begin{equation}
\label{a3.14}
\hhat S({\fg_-^*},V)\df\Sstar{\fg_-}\hotimes V
\end{equation}
of {\bf $V$-valued formal polynomials of $\fg_-^*$}.

Let $\eps_\0$, resp. $\eps_\1$, be an even, resp. odd, indeterminate.
Let the supercommutative superalgebra $A$ be defined as
$A=\K[\eps_\0,\eps_\1]/(\eps_\0^2)\otimes\Sstar{\fg_-}$.
Fix some basis $\{P_i\}_{i\in I}$ of $\fg_-$
and the dual 
basis 
$\{X^i\}_{i\in I}$ of $\fg_-^*$. Then, as follows from Corollary~\ref{SE*=K[X]}
in sect.~\ref{SE*}, one can identify the superalgebra $A$ with the superalgebra
$\K[[\eps_\0,\eps_\1,X]]/(\eps_\0^2)$ of formal polynomials in $\eps$ and $X$.

Let $\M$ be a homogeneous element of $\fg_-$, and $\eps=\eps_i$
if $\M\in{}_i\fg_-$. Let
\[
XP\df\sum_iX^i\otimes P_i\in
\Sstar{\fg_-}\hotimes \fg_-
\subset A\hotimes U(\fg). %
\]

Note that $XP$ is just the image of the identity map $\Id_{\fg_-}$
under the canonical 
morphism
$\Hom(\fg_-,\fg_-)\approx
\fg_-^*\hotimes\fg_-\subset\Sstar{\fg_-}\hotimes \fg_-$, so it does not depend
of the choice of the basis $\fam{P_i}{i\in I}$ and the dual basis $\fam{X^i}{i\in I}$.

In what follows we will often write just
$X^iP_i$ instead of $\sum_iX^i\otimes P_i$, etc., omitting $\otimes$
and assuming, (as is common practice in physical literature) that
summation is taken over repeated pairs of indices (subscript and superscrip), unless the
contrary is stated explicitly.

Consider an element
\[
\exp(\eps \M)\exp(XP)\in A\hotimes U(\fg).
\]

The exponents above exist, evidently, within the
superalgebra $A\hotimes U(\fg)$
.

Let $\A$ and $\B$ be morphisms of superspaces
\begin{equation}
\label{A}
\A:\fg\arr \Sstar{\fg_-}\hotimes\fg_-,\ \ \ \M\mapsto \A(\M),
\end{equation}

\begin{equation}
\label{B}
\B:\fg\arr \Sstar{\fg_-}\hotimes\fh\ \ \ \M\mapsto \B(\M).
\end{equation}

\begin{Th}
\label{p7}
The following properties {\rm(a)--(c)} of maps~$(\ref{A})$, $(\ref{B})$
are equivalent:

{\rm(a)} For any $\M\in\fg$, one has the equality of formal power series:
\begin{equation}
\label{a3.4}
e^{\eps \M}e^{XP}=e^{XP+\eps\A(\M)}e^{\eps\B(\M)}
\end{equation}

{\rm(b)} For any \fh-module $\rho:\fh\arr \End(V)$,
the representation \T \M takes the following form
in the realization~$(\ref{a3.14})$:
\begin{equation}
\label{3.14}
\T \M=\cD(\chA(\M))\otimes\One_V +\hhat\rho_{\Sstar{\fg_-}}(\chB(\M))\;,\ \ \M\in\fg,
\end{equation}
where \cD is defined by Eqs.~$(\ref{3.12})$--$(\ref{3.13})$ %
and $\hhat\rho_A$ 
for a linearly topologized commutative superalgebra $A$
is defined in Appendix~\ref{a1}
{\rm(see the definition after~(\ref{a1.18}))};
\chA and \chB means \A, resp.~\B twisted by main antiautomorphism
of the superalgebra \Sstar{\fg_-}, i.e. with $X$ replaced by $-X$ in the corresponding
formal series.

{\rm($\rm b'$)} The equalities~$(\ref{3.14})$ are valid for the adjoint representation of \fh;

{\rm(c)} The representation \I \M takes the following form
in the realization~(\ref{symspan}):
\begin{equation}
\label{3.14I}
\I \M=\cD_\infty(I\A(\M))\otimes\One_V +\rho_{ S(\fg_-)}(\B(\M))\;,\ \ \M\in\fg,
\end{equation}
where $\cD_\infty$ is defined by Eqs.~(\ref{3.12})--(\ref{3.13}) 
and $I:\Sstar{\fg_-}\hotimes\fg_-\arr\fg_-\hotimes\Sstar{\fg_-}$
is the canonical isomorphism of permutation;

{\rm($\rm c'$)} The equalities~$(\ref{3.14I})$ are valid for the adjoint representation of \fh.
\end{Th}

\begin{proof}

{\bf (a)$\implies$(${\rm}c'$)} Let
\begin{equation}
\label{varphi}
\A(\M)=\varphi^i(X,\M)P_i,\ \ \ \ \B(\M)=h(X,\M);
\end{equation}
so that
\begin{equation}
\label{4exp}
e^{\eps \M}e^{XP}=e^{(X+\eps\varphi(X,\M))P}e^{\eps h(X,\M)}\;.
\end{equation}

Expanding in powers of $X$ gives:
\[
e^{\varepsilon \M}e^{XP}=\sum_n{1\over n!}(1+\varepsilon \M)(XP)^n=
\sum_n{1\over n!}(1+\varepsilon \M)X^{\mu_1}P_{\mu_1}\cdots X^{\mu_n}P_{\mu_n}
\]

and
\[
e^{(X+\varepsilon\varphi(X,\M))P}e^{\varepsilon h(X,\M)}=
\sum_n{1\over n!}((X+\varepsilon\varphi(X,\M))P)^n(1+\varepsilon h(X,\M));
\]

i.e., linear on $\varepsilon$ terms in eq.~(\ref{a3.4}) will give
the following equality
\begin{equation}
\label{a3.7}
\sum_n{1\over n!}\M(XP)^n=
\sum_n{1\over n!}\left\{\sum_{k=1}^n(XP)^{k-1}\varphi P(XP)^{n-k}+
(XP)^nh(X,\M)\right\}.
\end{equation}

Further,
$(XP)^n=(XP)^{\circ n}$,
$\sum_{k=1}^n(XP)^{k-1}\varphi P(XP)^{n-k}=
n(\varphi P)\circ(XP)^{\circ(n-1)}$,
where $\circ$ means symmetric composition in \U{A\otimes\fg_-} defined
by eqn.(\ref{3.3}).

Let
\[
\renewcommand{\arraystretch}{1.4}
\begin{array}{ll}
\varphi^i(X,\M)=\sum_{\stackrel k{\nu_1\le\dots\le\nu_k}}
a^i_{\nu_1\dots\nu_k}X^{\nu_1}\cdots X^{\nu_k}\ \ &%
(a^i_{\nu_1\dots\nu_k}\in\K, |X^{\nu_1}|+\dots|X^{\nu_k}|+|\varepsilon|=|\varphi^i|),\\
h(X,\M)=\sum_{\stackrel k{\nu_1\le\dots\le\nu_k}}
h_{\nu_1\dots\nu_k}X^{\nu_1}\cdots X^{\nu_k}\ \ &%
(h_{\nu_1\dots\nu_k}\in\fh)
\end{array}
\]
be the expansions of $\varphi^i$ and $h$ into formal power series. 

Since
$|\varepsilon|+|\varphi^i|+|P_i|=0$
it follows that
\[
\varphi^iP_i=(-1)^{|\varphi^i||P_i|}P_i\varphi^i=(-1)^{(1+|\M|)|P_i|}P_i\varphi^i,
\]
and hence we obtain:
\begin{equation*}
\begin{split}
&\sum_n{1\over n!}\M(XP)^n=\\
&\sum_n{1\over(n-1)!}\sum
\left\{(-1)^{|\varphi^i|P_i|}P_i\circ
a^i_{\nu_1\dots\nu_k}X^{\nu_1}\cdots X^{\nu_k}(XP)^{\circ(n-1)}+
h_{\nu_1\dots\nu_k}\circ X^{\nu_1}\cdots X^{\nu_k}(XP)^{\circ n}\right\},
\end{split}
\end{equation*}
where
$h_{\nu_1\dots\nu_k}\circ\dots$
means that
$h_{\nu_1\dots\nu_k}$
must be placed on the right after all operations are performed.

Further, we have:
\begin{alphalabel}
\begin{equation}
\label{a3.8a}
{\partial\over\partial P_{\nu_1}}
\cdots
{\partial\over\partial P_{\nu_k}}(XP)^m=
{m!\over(m-k)!}X^{\nu_1}\cdots X^{\nu_k}(XP)^{m-k}
\end{equation}
or, equivalently:
\begin{equation}
\label{a3.8b}
X^{\nu_1}\cdots X^{\nu_k}(XP)^l=
{l!\over(k+l)!}
{\partial\over\partial P_{\nu_1}}
\cdots
{\partial\over\partial P_{\nu_k}}(XP)^{k+l}.
\end{equation}
\end{alphalabel}

Hence,
\begin{equation*}
\label{a3.9}
\begin{split}
&\sum_n{1\over n!}\M(XP)^{\circ n}=\\
&\sum_n{1\over(n-1)!}\sum
(-1)^{|\varphi^i||P_i|}P_i\circ a^i_{\nu_1\dots\nu_k}{(n-1)!\over(k+n-1)!}
{\partial\over\partial P_{\nu_1}}
\cdots
{\partial\over\partial P_{\nu_k}}(XP)^{k+n-1}\\
&+\sum_n{1\over n!}\sum
h_{\nu_1\dots\nu_k}\circ{n!\over(k+n)!}
{\partial\over\partial P_{\nu_1}}
\cdots
{\partial\over\partial P_{\nu_k}}(XP)^{k+n}\\
&=
\sum_n{1\over n!}\left\{(-1)^{|\varphi^i||P_i|}P_i\circ
\sum
a^i_{\nu_1\dots\nu_k}
{\partial\over\partial P_{\nu_1}}
\cdots
{\partial\over\partial P_{\nu_k}}\right.\\
&\left.+\sum
h_{\nu_1\dots\nu_k}\circ
{\partial\over\partial P_{\nu_1}}
\cdots
{\partial\over\partial P_{\nu_k}}\right\}\circ(XP)^{\circ n}\\
&=\sum_n{1\over n!}\left\{\sum_i(-1)^{|\varphi^i||P_i|}
P_i\varphi^i({\partial\over\partial P},\M)
+h({\partial\over\partial P},\M)\right\}\circ(XP)^{\circ n}.
\end{split}
\end{equation*}

Taking into account that
$(-1)^{|\varphi^i||P_i|}=(-1)^{(1+|\M|)|P_i|}$, we immediately get from the last equality
the following expression for the operator $L(\M)$ of left
shift by $\M$ in the space
\[
\U\fg\otimes_{\U\fh}\fh\approx\S\fg\otimes\fh
\]
as a differential operator (of infinite order, generally speaking):
\begin{equation}
\label{a3.10}
L(\M)=(-1)^{(1+|\M|)|P_i|}P_i\varphi^i({\partial\over\partial P},\M)+
h({\partial\over\partial P},\M).
\end{equation}
This proves ($\rm c'$).

The latter expression for a ``generic'' induced representation (i.e.,
induced from the adjoint representation of \fh) immediately implies
that for any 
representation $\rho$ of \fh in $V$
the operator $I_\fh^\fg(\M)$ of the action of $\M$ in the superspace of induced representation
\begin{equation}
\label{a3.11}
I(V)\equiv\U\fg\otimes_{\U\fh}V
\end{equation}
can be represented as
\begin{equation}
\label{a3.12}
\boxed{
I_\fh^\fg(\M)=(-1)^{(1+|\M|)|P_i|}P_i\varphi^i({\partial\over\partial P},\M)+
\rho(h({\partial\over\partial P},\M)),
}
\end{equation}

i.e., ($\rm c'$) $\implies$(c).

It is clear, that starting from (\ref{a3.10}) one can go backwards to get
the equality~(\ref{a3.7}), i.e., the implication ($\rm c'$)$\implies$(a) is also valid.

Now the expression~(\ref{3.14}) for the operators of coinduced
representation follows from that for the induced ones by duality,
if only the superspace $V$ is finite-dimensional,
due to an easily established (see~\cite{9}) canonical isomorphism
\begin{equation}
\label{caniso}
\I V^*\approx\T {V^*}.
\end{equation}

To extend the implication ($\rm a$)$\implies$(b) 
to the case where $V$ has infinite-dimension
(and is considered together with the discrete topology on it),
we have to find an expression for the right action of
an arbitrary homogeneous element of \fg on \U\fg in the
corresponding bases.

Let $\varphi\colon\S{\fg_-}\arr V$ be a linear map. Let
$\M\in\fg$ and $g'\in\U{\fg_-}\subset\U\fg$ be arbitrary.
Then 
the isomorphism~(\ref{PBWiso}) implies
that there exists a decomposition
\[
g'\M=\sum_i h_ig_i\ \ \ (h_i\in\U\fh,\ g_i,\in\U{\fg_-})
\]
such that
\[
(\T \M f)(g')=(-1)^{|\M|(|f|+|g'|)}\sum_i(-1)^{|h_i||f|}\rho(h_i)f(g_i).
\]

So we must find the expression for the right ``action'' of \fg
on $\U{\fg_-}\approx\S{\fg_-}$ first.

Inverting the equality~(\ref{a3.4}) one gets
\begin{equation*}
\label{a3.16}
e^{-XP}e^{-\eps \M}=e^{-\eps h(X,\M)}e^{-(X+\eps\varphi(X,\M))P}\;,
\end{equation*}
or, replacing $\eps\mapsto-\eps$, $X\mapsto-X$, we get
\begin{equation}
\label{a3.4'}
e^{XP}e^{\eps \M}=
e^{\eps h(-X,\M)}e^{(X+\eps\varphi(-X,\M))P}\;.\tag{\ref{a3.4}$'$}
\end{equation}

Now, the counterpart to eq.~(\ref{a3.7}) is:
\begin{equation}
\label{a3.7'}
\begin{split}
&\sum_n{1\over n!}(XP)^n\M=\\
&\sum_n{1\over n!}\left\{\sum_{k=1}^n(XP)^{k-1}\varphi^i(-X,\M)P_i(XP)^{n-k}+
(XP)^nh(-X,\M)\right\}\\
&=\sum_n{1\over n!}\left\{n\varphi^i(-X,\M)P_i(XP)^{\circ(n-1)}+
h(-X,\M)(XP)^{\circ n}\right\}.
\end{split}
\tag{\ref{a3.7}$'$}
\end{equation}

Expanding both $\varphi^i(-X,\M)$ and $h(-X,\M)$ in powers of $X$
gives the following result:
\begin{equation*}
\label{a3.9}
\begin{split}
&\sum_n{1\over n!}(XP)^{\circ n}\M=\dots\\
&=\sum_n{1\over n!}\left\{\sum_i(-1)^{(1+|\M|)|P_i|}
P_i\varphi^i(-{\partial\over\partial P},\M)
+h(-{\partial\over\partial P},\M)\right\}\circ(XP)^{\circ n}.
\end{split}
\end{equation*}

In other words, for the operator $R(\M)$ of right multiplication
of $\U{\fg_-}\approx\S{\fg_-}$ by $\M$ we get the following
expression:
\begin{equation}
\label{a3.10$'$}
R(\M)=(-1)^{(1+|\M|)|P_i|}P_i\varphi^i(-{\partial\over\partial P},\M)+
h(-{\partial\over\partial P},\M).
\tag{\ref{a3.10}$'$}
\end{equation}

Recall now that the action of an element
$f\otimes v\in\Sstar{\fg_-}\hotimes V$ considered as an element of
\L{\S{\fg_-}}V is defined by
\begin{equation}
\label{a3.15}
(f\otimes v)g=(-1)^{|v|+|g|}f(g)v
\end{equation}
(see Prop.~(\ref{WhotimesE})).
So, we have:
\begin{equation*}
\begin{split}
&\T\M(f\otimes v)(g')&=&(-1)^{|\M|(|f|+|v|+|g'|)}(f\otimes v)(g'g)\\
&&=&(-1)^{|\M|(|f|+|v|+|g'|)+|\M||g'|}(f\otimes v)(R(\M)g')\\
&&=&(-1)^{|\M|(|f|+|v|)}(f\otimes v)
((-1)^{(1+|\M|)|P_i|}P_i\varphi^i(-{\partial\over\partial P},\M)g'+h(-{\partial\over\partial P},\M)g')\\
&&=&(-1)^{|\M|(|f|+|v|)+|v|(|g'|+|\M|)}f
((-1)^{(1+|\M|)|P_i|}P_i\varphi^i(-{\partial\over\partial P},\M)g'+h(-{\partial\over\partial P},\M)g')v\\
&&=&(-1)^{|\M||f|+|v||g'|+|\M||f|}
((-1)^{(1+|\M|)|P_i|}P_i\varphi^i(-{\partial\over\partial P},\M)+h(-{\partial\over\partial P},\M))^*f(g')v\\
&&=&(-1)^{|v||g'|}(\varphi^i(-X,\M){\partial\over\partial X^i}+h(-X,\M))f(g')v\\
&&=&((-1)^{|X^i|}\varphi^i(-X,\M)){\partial\over\partial X^i}+h(-X,\M)))(f\otimes v)(g').
\end{split}
\end{equation*}

I.e., the action \T\M on elements of finite rank in \L{S(\fg_-)}V is given by
\begin{equation}
\label{a3.12*}
\T\M
=(-1)^{|X^i|}\varphi^i(-X,\M)){\partial\over\partial X^i}+h(-X,\M).
\end{equation}
But elements of finite rank are dense in \L{S(\fg_-)}V by heading (d) of Prop.~\ref{WhotimesE},
so that~(\ref{a3.12*}) is valid for any element of $\L{\fg_-}V=S(\fg_-)^*\hotimes V$.

The equality~(\ref{a3.12*}) is equivalent to the particular case of~(\ref{3.14})
by~(\ref{3.13expl}) and~(\ref{3.18}),
proving the implication (a)$\implies$($\rm b'$).
Other implications are proved similarly to already proved.
\end{proof}

The main question now is whether the maps $G_-$ and $H$ do
exist. And how they look like if they exist.

\section{The case where $\fg_-$ is a Lie superalgebra.}
\label{subcase}
In this section, let $\fg_-$ be a 
Lie superalgebra. 
In this case one can easily prove the existence of maps
(\ref{A}) and (\ref{B}) satisfying conditions of Theorem~\ref{p7}
as well as to find an explicit expression for \A and \B.

Let $\Pi_{\fg_-}$ and $\Pi_\fh$ be the projections of \fg onto $\fg_-$ and \fh, respectively:
\begin{equation}
\label{a3.2}
\Pi_{\fg_-}+\Pi_\fh=\One\;,\ \ \Pi_{\fg_-}^2=\Pi_{\fg_-}\;,\ \ %
\Pi_\fh^2=\Pi_\fh\;,\ \ \Pi_{\fg_-}\Pi_\fh=\Pi_\fh\Pi_{\fg_-}=0.
\end{equation}

Since $\eps^2=0$, we have for any $\M\in\fg$:
\begin{eqnarray*}
e^{\eps \M}e^{XP}&=&e^{XP}e^{-XP}e^{\eps \M}e^{XP}=
e^{XP}\exp(e^{-\ad XP}\eps \M)=\\
&&e^{XP}e^{\eps\Pi_{\fg_-}e^{-\ad XP}\M}
e^{\eps\Pi_\fh e^{-\ad XP}\M}\;.
\end{eqnarray*}

The same identity $\eps^2=0$ implies that (a well known fact)
\begin{equation}
\label{a3.3}
e^Ae^{\eps B}=\exp\left(A-\eps{\ad A\over e^{-\ad A}-1}B\right)
\end{equation}
for any $A$ and $B$ in any associative algebra where the
corresponding expressions make sense.

Hence, for any $\M\in\fg$
\begin{equation}
\label{a3.4xxx}
e^{\eps \M}e^{XP}=e^{(X+\eps\varphi(X,\M))P}e^{\eps h(X,\M)}\;,
\end{equation}
where
$\varphi^i(X,\M)\in\K[[X]]$ 
and
$h(X,\M)\in\K[[X]]\hotimes\fh$
are defined by
\begin{equation}
\label{a3.5}
\varphi^i(X,\M)P_i=
-{\ad XP\over e^{-\ad XP}-1}\One\otimes\Pi_{\fg_-}
e^{-\ad XP}\M\;;
\end{equation}
\begin{equation}
\label{a3.6}
h(X,\M)=
\One\otimes\Pi_\fh e^{-\ad XP}\M\;.
\end{equation}

In other words, in this case, the conditions of heading a) of Theorem~\ref{p7}
are satisfied with
\[
\A(\M)=\varphi^i(X,\M)P_i,\ \ \ \B(\M)=h(X,\M)
\]
so that representations~(\ref{3.14}), (\ref{3.14I}) for \T \M and \I \M are valid.

We note that if \fg is a \Z-graded Lie superalgebra
\begin{equation}
\label{3.20}
\fg=\bigoplus_{-l}^\infty\fg_i\;\ \ \ %
[\fg_i,\fg_j]\subset\fg_{i+j}\ \ i,j\in\Z
\end{equation}
and
\begin{equation}
\label{3.21}
\fh=\sum_{n=0}^\infty\fg_n\;,\ \ \ %
\fg_-=\sum_{n=-l}^{-1}\fg_n\;,
\end{equation}
such that the superalgebra $\fg_-$
is {\it finite\/}-dimensional,
then the infinite sums in the r.h.s. of
Eqs.~(\ref{a3.5}) and~(\ref{a3.6})
become finite 
and the differential operators~(\ref{3.14}) for $\M\in\fg_d$
have polynomial coefficients with degrees
of the corresponding polynomials not exceeding
$l+d$.

\section{The general case.}
It turns out that in general case,
where the subsuperspace $\fg_-$ need not be a subsuperalgebra of \fg,
the maps \A and \B
satisfying eq.~(\ref{a3.4}) exist as well. Here, in order to compactify
the formulas we use the shorthand expressions
$\Pi_{\fg_-}$ and $\Pi_\fh$ to denote $\One\otimes\Pi_{\fg_-}$
and $\One\otimes\Pi_\fh$,respectively.
\begin{Th}For any homogeneous $\M\in\fg$
the equality
\begin{equation}
\label{a3.4xxx'}
e^{\eps \M}e^{XP}=e^{(X+\eps\varphi(X,\M))P}e^{\eps h(X,\M)}\;
\end{equation}
is valid, where
\begin{equation}
\label{a3.5'}
\varphi^i(X,\M)P_i=
-\left.\left(\Pi_{\fg_-}{e^{-\ad XP}-1\over\ad XP}\right)\right|_{A\hotimes\fg_-}^{-1}\Pi_{\fg_-}
e^{-\ad XP}\M\;;
\end{equation}
\begin{equation}
\label{a3.6'}
h(X,\M)=
\left.\left(\Pi_\fh{\ad XP\over e^{-\ad XP}-1}\right)\right|_{A\hotimes\fh}^{-1}
\Pi_\fh{\ad XP\over e^{-\ad XP}-1}e^{-\ad XP}\M\;.
\end{equation}
\end{Th}

\begin{proof}
For the l.h.s. of~(\ref{a3.4xxx'}) we have due to the equality~(\ref{a3.3}):
\begin{equation}
\label{a3.noname}
e^{\eps \M}e^{XP}=e^{XP}e^{e^{-\ad XP}\eps \M}=
\exp\left(XP-{\ad XP\over e^{-\ad XP}-1}
e^{-\ad XP}\eps \M\right).
\end{equation}
On the other hand, the application of eq.~(\ref{a3.3})
to the r.h.s. of~(\ref{a3.4xxx'}) and taking into account that $\eps^2=0$, yields:
\begin{equation}
\label{general.1}
e^{(X+\eps\varphi(X,\M))P}e^{\eps h(X,\M)}=
\exp\left(XP+\eps\varphi(X,\M)P-{\ad XP\over e^{-\ad XP}-1}\eps h(X,\M)\right).
\end{equation}
Taking logarithm of right hand sides of~(\ref{a3.noname}) and~(\ref{general.1})
we see that the equation~(\ref{a3.5'}) is equivalent to the equation:
\begin{equation}
\label{general.2}
\varphi(X,\M)P-{\ad XP\over e^{-\ad XP}-1}h(X,\M)=
-{\ad XP\over e^{-\ad XP}-1}e^{-\ad XP}\M.
\end{equation}
Applying to both sides of the latter equation the operator
$\Pi_{\fg_-}{e^{-\ad XP}-1\over\ad XP}$ we obtain
\begin{equation}
\label{general.3}
\Pi_{\fg_-}{e^{-\ad XP}-1\over\ad XP}\varphi(X,\M)P=-\Pi_{\fg_-}e^{-\ad XP}\M.
\end{equation}
Clearly, the latter equation implies~(\ref{a3.5'}).

To obtain~(\ref{a3.6'}), apply $\Pi_\fh$ to both sides
of equation~(\ref{general.2}) to get
\begin{equation}
\label{general.6}
\Pi_\fh{\ad XP\over e^{-\ad XP}-1}h(X,\M)=
\Pi_\fh{\ad XP\over e^{-\ad XP}-1}e^{-\ad XP}\M.
\end{equation}
The latter equation implies evidently~(\ref{a3.6'}).
\end{proof}

\section{The action graph and ``path integral''.}
Expressions~(\ref{a3.5}) and~(\ref{a3.6}) are good enough
for manual calculations of \T \M and \I \M in the case of
\Z-graded Lie superalgebra of finite depth, considered
in section~\ref{subcase}: one expands the formal series
consisting of exponents (Bernoulli coefficients appear
at this stage); then one calculates multiple brackets
of the element $XP$ with \M. There is only a finite number
of non-zero  brackets for any fixed $\M$, so this process
terminates even in the case of infinite dimensional \fg
if only $\fg_-$ is finite-dimensional.

But to make all calculations manually is a tedious task
even for \fg of dimension a dozen or two. 

So it is desirable to convert expressions~(\ref{a3.5}) and~(\ref{a3.6})
into something potentially more digestible by computers.

A useful way to do this is by means of {\it action graphs\/} defined below.
Our definition will be a bit more general than we need in this paper,
so that it may be applied to 
local Lie superalgebras and Cartan-Tanaka-Shchepochkina prolongs (\cite{T}), (\cite{ALS}).

Given two superspaces $E$ and $H$,
any bilinear map
\begin{equation}
\label{action}
E\times H\stackrel\cdot\arr H
\end{equation}
is said to be an {\bf action}
of $E$ on $H$. Stress that neither $E$ nor $H$ are supposed to be equipped with
any other structure, except a superstructure.

Let $P=\fam{P_i}{i\in I}$ be a basis in $E$ and $h=\fam{h_j}{j\in J}$
be a basis in $H$.

In the bases $P$ and $h$, the action~(\ref{action})
can be expressed as follows:
\begin{equation}
\label{struct.const}
P_i\cdot h_j=c_{ij}^kh_k.
\end{equation}
The family $C=\fam{c_{ij}^k}{i,j\in I;k\in J}$ will be called
the {\bf structure constants} of the action~(\ref{action})
in bases $P$ and $h$.

Define the directed graph $G(P,h)$ as follows.
Its vertices are elements of the basis $h$;
for any vertices $h_s$ and $h_t$
and an element $P_i$ of the basis $P$ of $E$
there exists an edge from $h_s$ to $h_t$
labelled by $P_i$ 
iff $c_{is}^t\neq0$
.
Note, that there may be more than one  edge going from $h_s$ to $h_t$.

The graph 
$G(P,h)$
will be called the {\bf action graph}
of the action~(\ref{action}) with respect to the bases $P$ and~$h$.

Before proceeding further we are to give some (ad hoc) definitions.

For any directed graph $G$ denote $\cP(G)$ the set of all directed paths of $G$.

$\cP(G)$, equipped with partial operation of composition of paths,
clearly may be considered as a category with the set of objects coinciding
with the set of vertices of $G$.

The set of all paths beginning at vertex $v$, resp. ending at vertex $v'$,
resp. beginning at $v$ and anding at $v'$ will be denoted
$\cP_v(G)$, resp. $\cP^{v'}(G)$, resp. $\cP_v^{v'}(G)$.

Let $A$ be a supercommutative superalgebra with unit element. A map
\begin{equation}
\label{muA}
\mu:\cP(G)\arr A
\end{equation}
such that
$\mu(p\circ p')=\mu(p)\mu(p')$ and $\mu(1_v\circ p)=\mu(p)$
(i.e. $\mu$ is a functor into multiplicative monoid of $A$)
will be called a {\bf multiplicative $A$-valued measure} on $\cP(G)$.
If, besides, $V$ is an $A$-supermodule, then the map
\begin{equation}
\label{muV}
\mu_V:\cP(G)\arr V
\end{equation}
such that
$\mu_V(p\circ p')=\mu(p)\mu_V(p')$
will be called a {\bf $\mu$-compatible multiplicative $V$-valued measure}.

Clearly, a multiplicative $A$-valued measure on $\cP(G)$
is uniquely determined by its values on edges of the graph $G$
(i.e., paths of length 1). Moreover, any map
$
\text{Edges}(G)\arr A
$ 
can be extended to a multiplicative $A$-valued measure on $\cP(G)$.

To define path integrals over a graph, suppose first
that the graph $G$ is such that, for any vertex $v$ of $G$,
there exists only {\sl finite} number of directed paths starting at $v$.
This is the case of the action graph of the action $\fg_-\times\fg\arr\fg$,
where \fg is a \Z-graded Lie superalgebra and $\fg_-=\oplus_{i<0}\fg_i$
is finite-dimensional,
considered at the end of preceeding section.

Then, for any function $c:\cP(G)\arr A$, there exists
\begin{equation}
\label{int}
\int_{v}c\mu_V\df\sum_{p\in\cP_v(G)}c(p)\mu_V(p).
\end{equation}

The sum in the r.h.s. of definition~(\ref{int})
may have sense (for particular choice of both $\mu_V$ and $c$)
even if it is infinite in the case where
$A$ is a {\it topological} superalgebra
and $V$ is a {\it topological} supermodule.

We will not give here the most general definitions,
restricting ourselves
to the case of action graph $G(P,h)$ with a special kind of
multiplicative measure defined on it.
Namely, let $A=\Sstar E$ and $V=\Sstar E\hotimes H$.
Let \fam{X_i}{i\in I} be the (topological) basis in $E^*$ dual
to the basis $P$ in $E$.

Let $\mu$ be the multiplicative $A$-valued measure on
$\cP(G(P,h))$ which extends the map
\begin{equation}
\label{muA'}
\mu_0:\text{Edges}(G(P,h))\arr A\ \ \ %
((h_s\stackrel{P_i}\arr h_t)\mapsto -c_{is}^tX_i).
\end{equation}
Let $\mu_V$ be the multiplicative $V$-valued measure compatible
with $\mu$ and defined as follows:
\begin{equation}
\label{muV'}
\mu_V(p)\df\mu(p)T(p),
\end{equation}
where 
$T(p)$ denotes the end (i.e., target) vertex of the path $p$.

\begin{Prop}
For any action $E\times H\arr H$ and any bases $P$ in $E$ and $h$ in $H$,
the sum in the r.h.s. of eq.~$(\ref{int})$ exists
for the measure~$\mu_V$ defined by eq.~$(\ref{muV'})$
and any function
\[
c:\cP(G(P,h))\arr\K\harr\Sstar E\hotimes H.
\]
\end{Prop}

\begin{proof}
This trivially follows from the definition of
formal series module $S(E^*)\hotimes H$ and
the fact that any monomial of degree $n$
may occur in the sum~(\ref{int}) only finite number of times (not more than $n!$).
\end{proof}

The next theorem states that the formal series $\varphi^i(X,\M)$
and $h(X,\M)$ defined by eqs.~(\ref{a3.5})--(\ref{a3.6})
giving the expressions for (co)induced modules in case
where $\fg_-$ is a subsuperalgebra of \fg can be calculated as
path integral on the action graph with an appropriate function $c(p)$.

For nonnegative integers $k$ and $n$ such that $k\leq n$, define
the rational number $c(k,n)$ as follows:
\begin{equation}
\label{c(k,n)}
c(k,n)=\sum_{k\leq i\leq n}{b_{n-i}\over i!(n-i)!},
\end{equation}
where $b_i$ is the $i$-th Bernoulli number.

For a path $p\in\cP(G(P,h))$, let $k(p)$ be the length of the longest
subpath, whose source vertex coincides with that of $p$,
whereas the target vertex belongs to \fh (recall that vertices
of $\cP(G(P,h))$ are basis elements of \fg).
Define now the map $K:\cP(G(P,h))\arr\K$ as follows:
\begin{equation}
\label{Kdef}
K(p)\df-c(k(p),\text{length}(p)),
\end{equation}
where $\text{length}(p)$ denotes the number of edges in the path $p$.

\begin{Th}
If $\fg_-$ is a Lie subsuperalgebra of \fg,
then, for any basis element $M\in P\cup h$,
the 
formal series
$\A(M)$ and $\B(M)$
\DL{SKAZHI KTO TAKAYA $M$ NE V OTVETE TRUDYASHEMU, A V TEXTE!}
\VM{Dyk ono i skazano vrode b: ``basis element'' :))}
satisfying condition~$(\ref{a3.4})$ are determined as follows:
\begin{equation}
\label{A&B}
\A(M)=\Pi_{\fg_-}\int_M K\mu_V,\ \ \ \ \B(M)=\Pi_\fh\int_M K\mu_V,
\end{equation}
where $\mu_V$ is defined by eq.~$(\ref{muV'})$
and $K$ is defined by eq.~$(\ref{Kdef})$, and the projections
$\Pi_{\fg_-}$ and $\Pi_\fh$ are  defined by
eq.~$(\ref{a3.2})$.
In other words, $\A(M)$ is the sum over all paths starting at $M$
and ending in $\fg_-$, whereas $\B(M)$ is the sum over all paths starting at $M$
and ending in \fh.
\end{Th}

\begin{proof}
Meditate over eqs.~(\ref{a3.5}) and~(\ref{a3.6}) until being enlightened :).
\end{proof}

\section{The program {\tt coindrep}.}
The formulae~(\ref{A&B}) were converted into algorithm
written in ANSI C++. The present beta version of the program
{\tt coindrep} which uses this algorithm calculates
generators of series $\mathfrak{gl}(n)$, $A_n=\mathfrak{sl}(n+1)$,
$D_n=\mathfrak{so}(2n)$ and $E_n$ only, for the case
where \fh in the expansion $\fg=\fg_-\oplus\fh$ is a maximal solvable
subalgebra of \fg. The possibility to make calculations for other series of
simple Lie algebras as well as Lie superalgebras will be added to the program
as well as the possibility to make calculations for arbitrary
``good'' expansion $\fg=\fg_-\oplus\fh$ defined by user.

The program output (statistics and expressions for generators
of coinduced module) is in plain \TeX~format.

For the latest version of the windows executable
of the program together with its source code, see\linebreak
{\bf theo.inrne.bas.bg/\~{}vmolot/coindrep.zip}.

\subsection{Some statistics.}
Calculations for the case $\fg=\mathfrak{gl}(15)$ take about half an hour
on Sun Solaris 9 system. The program is not optimized yet, but even
for optimized future version this time is expected to be proportional
to the number of paths in the action graph, so that, roughly, will grow
exponentially as the depths of canonical \Z-grading grows.

The longest number of paths in path integral for the action graph of
$\fg=\mathfrak{gl}(15)$ is 3052080, 
whereas the number of monomials
after reduction is exactly 8K=8192 for all 15 generators from $\fg_+$.
The maximal degree of these monomials is 15.

Statistics for the Lie algebra $E_6$: the longest number of paths is 73179,
the maximal number of monomials after reduction is 1906, the maximal degree
of monomials of generators from $\fg_1$ is 12. All calculations take less than 10 seconds.

In contrast, calculations of all generators for Lie algebra $E_7$ are estimated to take
about 10--20 hours. I have no guts to run the program to the end, restricting myself
to the calculations of the first few ``easy'' generators of maximal degree 16, which
took about half an hour. I hope that after optimization of the program
(replacing, where appropriate, STL\footnote{Standard Template Library of C++}
vectors by data structures based on trees,
like maps and multisets) the run time can be reduced essentially.
But I suspect that, even after eventual optimization, the calculations for
the Lie algebra $E_8$ will occupy many years of present supercomputer's time.
This is due to the fact that this time
is expected to be proportional to the exponent of
the depth of the main \Z-grading for $E_8$ which is 29.

\appendix
\section{Concise Dictionary on Lie superalgebras.}
\label{a1}

Here we have placed together some basic definitions
and constructions. For a detailed presentation of
``the linear algebra in superspace'', see
Chapter~$1$ of Ref.~\cite{7} and that of Lie superalgebras can be obtained
from reading Ref.~\cite{8}.

It
is supposed here that all vector spaces are over
some field \K of characteristic 0.
$\Zz=\Z/2\Z=\{\0,\1\}$,
$(-1)^\0=1$ and $(-1)^\1=-1$.

Let $V$ be a vector space. A {\bf superstructure}
on $V$ is defined by any of the three equivalent
structures:
\begin{list}{}{}
\item[a)] the linear operator $P\in\End_\K(V)$, called the {\bf parity operator}, such that
\begin{equation}
\label{a1.1}
P^2=\One_V\;;
\end{equation}
\item[b)] two linear operators ${}_i\Pi\in\End_\K(V)$
($i=\0,\1$) such that
\begin{equation}
\label{a1.2}
{}_i\Pi^2={}_i\Pi\;,\ \ {}_\0\Pi{}_\1\Pi={}_\1\Pi{}_\0\Pi=0\ \ %
\mbox{\rm and}\ \ {}_\0\Pi+{}_\1\Pi=\One_V\;;
\end{equation}
\item[c)] the decomposition
\begin{equation}
\label{a1.3}
V={}_\0V\oplus{}_\1V
\end{equation}
\end{list}
related to each other as follows:
\begin{equation}
\label{a1.4}
{}_i\Pi={\One+(-1)^iP\over2}\;,\ \ P={}_\0\Pi-{}_\1\Pi\;,\ \ %
{}_iV={}_i\Pi V\ \ (i=\0,\1).
\end{equation}

A {\bf superspace} is a linear space $V$ together with some
superstructure on it. For any $x\in V$,
the element ${}_\0x={}_\0\Pi x$ is said to be the {\bf even} component of $x$
and ${}_\1x={}_\1\Pi x$ is the {\bf odd} component.
If two or more superspaces are considered simultaneously,
the operators $P$ and ${}_i\Pi$ are sometimes supplied
with indices.

Here are standard constructions of new
superspaces from given ones. Given superspaces
$V$ and $W$, the following superstructures are naturally defined
on $V\otimes W$ and ${\rm L}_\K(V,W)$ (linear maps from
$V$ to $W$), respectively:
\begin{alphalabel}
\begin{eqnarray}
\label{a1.5}
P(v\otimes w)\df Pv\otimes Pw\ \ \ (v\in V,w\in W),\\
P(F)\df P_W\circ F\circ P_V\ \ \ (F\in {\rm L}_\K(V,W)).
\end{eqnarray}
\end{alphalabel}

In particular, the space $V^*={\rm L}_\K(V,k)$ is a superspace
if $V$ is a superspace.

A {\bf superalgebra} is a superspace $A$ together with
a bilinear map (composition):
\begin{equation}
\label{a1.6}
A\times A\arr A\ \ \ (x,y)\mapsto x\mult y\ \ \ (x,y\in A)
\end{equation}
such that
\begin{equation}
\label{a1.7}
P(x\mult y)=Px\mult Py\ \ \ (x,y\in A).
\end{equation}
Associativity and the identity elements in superalgebras are defined as
in ordinary algebras. A superalgebra $A$ is said to be {\bf supercommutative}
if
\begin{equation}
\label{a1.8}
xy=\sum_{i,j\in\Zz}(-1)^{ij}{}_iy\;{}_jx\ \ \ (x,y\in A).
\end{equation}

\noindent{\bf Examples.} If $V$ is a superspace, then
the tensor algebra
\[
T(V)=\bigoplus_0^\infty T^n(V)
\]
is an associative superalgebra with unity and
{\bf symmetric superalgebra}
\[
S(V)=\bigoplus_0^\infty S^n(V)
\]
of $V$ defined as the quotient algebra of $T(V)$
modulo relations
\[
x\otimes y-\sum_{i,j\in\Zz}(-1)^{ij}{}_jy\otimes{}_ix
\]
is a supercommutative associative superalgebra with identity.

\iffalse
Another example of supercommutative superalgebras is
the formal symmetric algebra
\[
\hhat S(V)=\prod_0^\infty S^n(V)
\]
of $V$. Superalgebras of this sort arise
in Sec.~\ref{formal}.\DL{?? vykin' poslednyuyu storku i voprosa ne budet...}
\else
The superspace $S(V)^*$ dual to the superalgebra $S(V)$ is a supercommutative superalgebra
as well (see section~\ref{SE*}). Moreover, it is a {\it topological} superalgebra,
if equipped with weak topology, turning $S(V)^*$ into linearly compact space.
\fi

A superalgebra \fg is a {\bf Lie superalgebra} 
if the composition in \fg, denoted usually as $[\cdot ,\cdot ]$,
satisfies the following relations:
\begin{equation}
\label{a1.9}
[x,y]=-\sum_{i,j\in\Zz}(-1)^{ij}[{}_iy,{}_jx]\ \ \ (x,y\in\fg)
\end{equation}
and
\begin{equation}
\label{a1.10}
(-1)^{ij}[{}_ix,[{}_jy,{}_kz]]+
(-1)^{jk}[{}_jy,[{}_kz,{}_ix]]+
(-1)^{ki}[{}_kz,[{}_ix,{}_jy]]=0\ \ (x,y,z\in\fg\ i,j,k\in\Zz).
\end{equation}

The {\bf tensor product} of the superalgebras $A$ and $B$
is the superspace $A\otimes B$ with composition defined
as follows:
\begin{equation}
\label{a1.11}
(a\otimes b)(a'\otimes b')=
\sum_{i,j\in\Zz}(-1)^{ij}a\;{}_ia'\otimes{}_jbb'\ \ %
(a,a'\in A\ b,b'\in B).
\end{equation}

The tensor product of two associative (supercommutative)
superalgebras is an associative (supercommutative) superalgebra.
The tensor product of a supercommutative superalgebra and
of Lie superalgebra is Lie superalgebra.

\noindent{\bf Examples.} 1) If $V$ is a superspace, then
we have the isomorphism of supercommutative superalgebras:
\begin{equation}
\label{a1.12}
S(V)\approx S({}_\0V)\otimes\La({}_\1V),
\end{equation}
where $\La({}_\1V)$ is the exterior algebra
of the vector space ${}_\1V$ supplied with
a natural superstructure---the extension
of the superstructure
$Px=-x$ of ${}_\1V$ by means of relations~(\ref{a1.7})

2) If \fg is a Lie superalgebra and $V$ is a superspace,
then $S(V)^*\otimes\fg$ (as well as its completion $S(V)^*\hotimes\fg$)
is a topological Lie superalgebra (see sect.~\ref{lpco}).

A {\bf linear representation} $\rho$ of a Lie superalgebra \fg
in a superspace $V$ is a linear map
\[
\rho\colon\fg\arr\End_\K(V)\df{\rm L}_\K(V,V)
\]
such that
\begin{equation}
\label{a1.13}
\rho([g_1,g_2])=\rho(g_1)\rho(g_2)-
\sum_{i,j\in\Zz}(-1)^{ij}\rho({}_ig_2)\rho({}_jg_1)\ \ (g_i\in\fg)
\end{equation}
and
\begin{equation}
\label{a1.14}
P(\rho(g))=\rho(P(g))\ \ \ (g\in\fg).
\end{equation}

The representatation $\rho^*$ in the superspace $V^*$
defined as
\begin{equation}
\label{a1.15}
\rho^*(g)=-(\rho(g))^*
\end{equation}
is said to be {\bf contragradient} to $\rho$. Here
the operator
$A^*\colon W^*\arr V^*$
{\bf dual} to the operator
$A\colon V\arr W$
is defined as follows
($\la\cdot ,\cdot \ra$ is the ``scalar product''
in $V^*\times V$ or in $W^*\times W$):
\begin{equation}
\label{a1.16}
\la A^*w^*,v\ra=
\sum_{i,j\in\Zz}(-1)^{ij}\la{}_iw^*,{}_jAv\ra\ \ \ (v\in V, w^*\in W^*).
\end{equation}
If $B\colon W\arr Y$ is another operator, then
\begin{equation}
\label{a1.17}
(B\circ A)^*=\sum_{i,j\in\Zz}(-1)^{ij}{}_iB^*\circ{}_jA^*.
\end{equation}

Given a representatation $\rho$ of a Lie superalgebra \fg,
we can define another representatation,
$\rho'(g)\df\rho(Pg)$,
where $P$ is the parity operator of the Lie superalgebra \fg.

Given a representatation $\rho$ of a Lie superalgebra \fg
in a superspace $V$ and any supercommutative associative
superalgebra $A$, define the representatation $\rho_A$
of the Lie superalgebra $A\otimes \fg$ in the superspace
$A\otimes V$ as follows:
\begin{equation}
\label{a1.18}
\rho_A(a\otimes g)(a'\otimes v)=
\sum_{i,j\in\Zz}(-1)^{ij}a\;{}_ia'\otimes{}_jgv\ \ %
(a,a'\in A,\ g\in G,\ v\in V).
\end{equation}

If, moreover, $A$ is a linearly topologized superalgebra,
then the representation~(\ref{a1.18})extends by continuity
to the representation $\hhat\rho_A$ of the linearly topologized
Lie superalgebra $A\hotimes \fg$
in the superspace $A\hotimes V$.

Given two representatations $\rho$ and $\rho'$ of
a Lie superalgebra \fg in superspaces $V$ and $V'$,
respectively, we say that an operator $A\colon V\arr V'$ is
an {\bf intertwining operator} between representations
$\rho$ and $\rho'$ if
\begin{equation}
\label{a1.19}
\rho'(g)\circ A=\sum_{i,j\in\Zz}(-1)^{ij}{}_iA\circ\rho({}_jg).
\end{equation}


\begin{thebibliography}{aaa}
\iffalse
\bibitem[1]{1} V. Molotkov, ICTP preprint IC/80/104  (1980);
scanned version available at {\bf http://ccdb3fs.kek.jp/cgi-bin/img\_index?8012175}.
\bibitem[2]{2} I. L. Kantor,Formulas for infinitesimal operators
of the Lie algebra of a homogeneous space, Trudy Sem. Vektor Tensor Anal. {\bf 17}, 243--249 (1974).
\bibitem[3]{3}
I. L. Kantor, On a vector field formula for the Lie algebra of a
homogeneous space.  J. Algebra 235 (2001), no.  2, 766--782
\bibitem[4]{7}
D. A. Leites, Introduction to the supermanifold theory.  Russian Math.
Surveys, v.  35, 1980, no.  1, 3--53.
\bibitem[5]{8}
V. G. Kac, Lie superalgebras, Adv. Math {\bf 26}, No.1, 8--96 (1977).
\bibitem[6]{9}
R. J. Blattner, Induced and produced representations of Lie algebras.
Trans. Amer. Math. Soc.  144  (1969), 457--474.
\bibitem[7]{10}
V. G. Kac,  Representations of classical Lie superalgebras.  Differential
geometrical methods in mathematical physics, II (Proc. Conf., Univ. Bonn, Bonn, 1977),
pp. 597--626, Lecture Notes in Math., 676, Springer, Berlin, 1978.
\bibitem[8]{11} L.Corwin, Y.Ne'eman, S.Sternberg, Rev.Mod.Phys. {\bf 43}, 573 (1975).
\bibitem[9]{12}
J. N. Bernstein, D. A. Leites, Irreducible representations
of the Lie superalgebra of vector fields and invariant differential
operators. Serdika, Bulgarian Math.  J., v.
7, 1981, 320--334 (in Russian, English translation in Selecta
Math.  Sov., v.  1, 1981, no.  2, 143--160
\bibitem[10]{16}
N. Bourbaki, {\sl Commutative algebra. Chapters $1-7$}. Translated from the French.
Reprint of the 1989 English translation.
Elements of Mathematics (Berlin). Springer-Verlag, Berlin, 1998. xxiv+625 pp.
\bibitem[11]{DG}
W. A. De Graaf {\sl Lie Algebras: Theory and Algorithms \/}, North-Holland Math. Library,
56. North-Holland Publishing Co., Amsterdam, 2000. xii+393 pp.
\bibitem[12]{T}
Tanaka N., On infinitesimal automorphisms of Siegel domains, J. Math.
Soc. Japan 22, 1970, 180--212;
\bibitem[13]{ALS} 
Alekseevsky D., Leites D., Shchepochkina I., New examples of simple
Lie superalgebras of vector fields.  C.r. Acad.  Bulg.  Sci., v.  34,
N 9, 1980, 1187--1190 (in Russian)
\else
\bibitem[1]{ALS} 
Alekseevsky D., Leites D., Shchepochkina I., New examples of simple
Lie superalgebras of vector fields.  C.r. Acad.  Bulg.  Sci., v.  34,
N 9, 1980, 1187--1190 (in Russian); for more details, see Shchepochkina I.,
Five exceptional simple Lie superalgebras of
vector fields and their fourteen regradings. Representation Theory
(electronic journal of AMS), v. 3, 1999, 3 (1999), 373--415
\bibitem[2]{12}
J. N. Bernstein, D. A. Leites, Irreducible representations
of the Lie superalgebra of vector fields and invariant differential
operators. Serdika, Bulgarian Math. J., v.
7, 1981, 320--334 (in Russian, English translation in Selecta
Math.  Sov., v.  1, 1981, no.  2, 143--160
\bibitem[3]{9}
R. J. Blattner, Induced and produced representations of Lie algebras.
Trans. Amer. Math. Soc.  144  (1969), 457--474.
\bibitem[4]{B1}
N.~Bourbaki, {\sl Topologie g\'en\'erale}, Hermann, Paris (1960).
\bibitem[5]{16}
N. Bourbaki, {\sl Commutative algebra. Chapters $1-7$}. Translated from the French.
Reprint of the 1989 English translation.
Elements of Mathematics (Berlin). Springer-Verlag, Berlin, 1998. xxiv+625 pp.
\bibitem[6]{S}
Sophie Chemla, Propri\'et\'es de dualit\'e dans les repr\'esentations coinduites
des superalg\`ebres de Lie. Ann.~Inst.~Fourier, Grenoble, {\bf44}, 4 (1994), 1067--1090.
\bibitem[7]{11} L.Corwin, Y.Ne'eman, S.Sternberg, Graded Lie algebras in mathematics and
physics (Bose-Fermi symmetry).
Rev. Modern Phys. 47 (1975), 573--603
\bibitem[8]{FrJa} A. Fr\"olicher, W. Jarchow, Zur Dualit\"atstheorie kompakt
erzeugter und lokalkonvexer Vektorr\"aume. Comm.~Math.~Helv., v.47, 289--310 (1972).
\bibitem[9]{DG}
W. A. De Graaf {\sl Lie Algebras: Theory and Algorithms \/}, North-Holland Math. Library,
56. North-Holland Publishing Co., Amsterdam, 2000. xii+393 pp.
\bibitem[10]{8}
V. G. Kac, Lie superalgebras, Adv. Math {\bf 26}, No.1, 8--96 (1977).
\bibitem[11]{10}
V. G. Kac,  Representations of classical Lie superalgebras.
Differential geometrical methods in mathematical physics, II (Proc. Conf., Univ. Bonn, Bonn, 1977),
pp. 597--626, Lecture Notes in Math., 676, Springer, Berlin, 1978.
\bibitem[12]{2} I. L. Kantor, {\sl Formulas for infinitesimal operators
of the Lie algebra of a homogeneous space\/}, Trudy Sem. Vektor Tensor Anal.{\bf 17}, 243-249 (1974).
\bibitem[13]{3}
I. L. Kantor, On a vector field formula for the Lie algebra of a
homogeneous space.  J. Algebra 235 (2001), no.  2, 766--782
\bibitem[14]{Ko}G. K\"othe, Topologische lineare R\"aume, I,
Die Grundlehren der Math. Wiss., bd.107, Springer-Verlag (1960).
\bibitem[15]{7}D. A. Leites, Introduction to the supermanifold theory.  Russian Math.
Surveys, v.  35, 1980, no.  1, 3--53.
\bibitem[16]{1} V. Molotkov, Equivalence between representations of conformal superalgebra
induced from different subgroups: invariant subspaces.
ICTP preprint IC/80/104  (1980);
scanned version available at {\bf http://ccdb3fs.kek.jp/cgi-bin/img\_index?8012175}.
\bibitem[17]{M2} V. Molotkov, Convenient categories of linearly topologized vector spaces
(in preparation);
\bibitem[18]{M3} V. Molotkov, Infinite-Dimensional Supermanifolds.
Trieste report IC/84/183 (1984);
\bibitem[18]{Po} H.~Porta, Compactly determined locally convex topologies.
Math.~Ann., 196, 91--100 (1972).
\bibitem[19]{T}
Tanaka N., 
On the equivalence problems
associated with simple graded Lie algebras. Hokkaido Math. J., 8
(1979), no. 1, 23--84
\bibitem[20]{Sh}
Schaefer H. {\sl Topological vector spaces.}
Third printing corrected. Graduate Texts in Mathematics, Vol. 3.
Springer-Verlag, New York-Berlin, 1971. xi+294 pp.
\fi
\end{thebibliography}
\end{document}